\documentclass{amsart}
\usepackage{myAMSpackages}
\usepackage{Universal}
\usepackage{MyAMSEnvironments}
\usepackage{Hladky}

\theoremstyle{plain}

\theoremstyle{plain}
\newtheorem{cnds}{Condition}
\renewcommand{\thecnds}{\Alph{cnds}}

\begin{document}

\author{Robert K. Hladky}
\address{University of Rochester, Rochester, NY 14627}
\email{hladky@math.rochester.edu}

\title{The \db-Neumann problem on noncharacteristic domains}

\footnotetext{\textit{Math Subject Classifications.}
			32V05, 35H20.}



\begin{abstract}
We study the $\bar{\partial}_b$-Neumann problem for domains $\Omega$  contained in a strictly pseudoconvex manifold $M^{2n+1}$ whose boundaries are noncharacteristic and have defining functions depending solely on the real and imaginary parts of a single CR function $w$. When the Kohn Laplacian is a priori known to have closed range in $L^2$, we prove sharp regularity and estimates for solutions. We establish a condition on the boundary $\partial \Omega$ which is sufficient for $\boxb$ to be Fredholm on $L^2_{(0,q)}(\Omega)$ and show that this condition always holds when $M$ is embedded as a hypersurface in $\cn{n+1}$.  We present examples where the inhomogenous $\bar{\partial}_b$ equation can always be solved in $C^\infty(\overline{\Omega})$ on $(p,q)$-forms with $1 \leq q \leq n-2$.
 \end{abstract}

\maketitle

\section{Introduction}\setS{ID}

In this paper we shall explore the boundary regularity of solutions to the $\db$-Neumann problem on compact domains inside strictly pseudoconvex pseudohermitian manifolds. We shall require that our domain $\Omega$ be noncharacteristic and satisfy the following condition

\renewcommand{\thecnds}{(A\arabic{cnds})}
\begin{cnds}
$\Omega$ posses a defining function $\varrho$ depending upon the real and imaginary parts of a particular CR function $w$.
\end{cnds}

The defining function part of this condition is typical for most discussions of solvability for either the Kohn Laplacian $\boxb$ or the $\db$-complex on domains in CR manifolds. The study of the $\db$-Neumann problem started with Kuranishi  (\cite{Kuranishi1}, \cite{Kuranishi2} and \cite{Kuranishi3}), who established existence for a weighted Neumann problem on small balls, as part of his study of the embeddability of strictly pseudoconvex CR structures. More recently, Shaw has established unweighted $L^2$-existence results for small sets of CR manifolds embedded in $\cn{n}$ whose defining function satisfies  (A1) and is convex in $w$ , see \cite{Shaw:B} or \cite{Shaw:P}. With the additional simplifying condition that the boundary has no characteristic points, Diaz has refined the techniques first employed by Kuranishi. In \cite{Diaz} he established that under the assumption of particular pointwise curvature bounds, $L^2$ solutions exist with exact Sobolev regularity for a problem closely related to the $\db$-Neumann problem. His solutions are only guaranteed to meet the second Neumann boundary condition. Exact regularity refers to estimates of the type $\|\varphi\|_{H^k} \leqc \|\boxb \varphi\|_{H^k}$. Diaz was interested in the tangential Cauchy-Riemann equations and his results are sufficient to show the existence of smooth solutions. However in general, the solutions exhibit a loss of Sobolev regularity. 

The analysis of the $\db$-Neumann problem is intricate as the operator \boxb is only subelliptic rather than elliptic. In addition the boundary conditions for the Neumann problem are non-coercive in the sense that the interior subelliptic estimates do not extend to the boundary of the domain $\Omega$. The presence of characteristic points on the boundary also complicates $L^2$ arguments enormously; the dimension of the horizontal space tangent to the boundary jumps. Both Kuranishi's and Diaz's argument for regularity involved the use of a subelliptic gain in directions tangent to the foliation by level sets of $w$. 

Improvements on these results were obtained by the first author in some special cases. Existence and sharp regularity were proved in \cite{hladky1} assuming (A1) and that the level sets of $w$ were all CR diffeomorphic to the same compact normal pseudohermitian manifold. This work was extended in \cite{hladky2} to include the homogeneous unit ball in the Heisenberg group. Additionally, some new negative results were obtained: the Kohn Laplacian can have infinite dimensional kernel, its partial inverse is non compact, and the Kohn Laplacian need not hypoelliptic.

On compact manifolds however, the Kohn Laplacian is well understood. In \cite{Folland:H} Folland and Stein introduced a new class of function spaces $S^k$ and proved sharp estimates for the Kohn Laplacian in terms of these. The first authors previous work involved decomposing the operator \boxb into pieces tangential and transverse to the foliation. Then it was possible to use global estimates on the compact leaves of the foliation and local elliptic estimates in the transverse directions to obtain sharp regularity and existence results.

In this paper, we shall extend the results of the first author's previous work to more general spaces. Namely, we shall establish sharp regularity for any noncharacteristic domain satisfying (A1) on which $\boxb$ has closed range on $L^2$. In particular, this means that our results will apply to the cases studied by Shaw \cite{Shaw:P} or Diaz \cite{Diaz}. The lack of uniformity of the foliating leaves is the main issue with the generalization. Especially note, that there may be jumps in the cohomology of the leaves which complicates estimates enormously. Additionally, this lack of uniformity means that we cannot decompose $\boxb$ into elliptic operators on hyperbolic space as in \cite{hladky1} and \cite{hladky2}, so instead we employ the technique of elliptic regularization. The key step in establishing the required \text{a priori} estimates is to adapt the interpolation techniques of \cite{Shaw:Int} to work simultaneously on all leaves.

 For technical reasons, the regularity results for the Kohn Laplacian will require an extra condition

 
\begin{cnds}
Near the boundary $\partial \Omega$,  $H_{\bt{w}} w$ is constant on leaves of the foliation by $w$ 
\end{cnds}

Here $H_{\bt{w}}$ is the pseudoHamiltonian of $w$ with respect to the pseudohermitian form $\eta$, defined in \rfS{NL}. This condition is non-generic, but if $\Omega$ has noncharacteristic boundary, there is a guaranteed pseudohermitian form on $M$ such that  (A2) holds.

Our main theorem is

\begin{thms}\label{A}
Let $\Omega$ be a smoothly bounded domain in a strictly pseudoconvex pseudohermitian manifold $(M,\eta)$ of dimension $2n+1$ with $n \geq 3$ such that $\Omega$ has noncharacteristic boundary and $(\Omega,w,\eta)$ satisfy (A1) and (A2).  Let  $1 \leq q \leq n-2$. If the Kohn Laplacian on $\Omega$ can be shown to have closed range as an unbounded operator on $L_{(0,q)} ^2 (\Omega)$ then following holds:

\hfill

\noindent For any  $(0,q)$-form $f$, there exists a unique solution $u \perp \Ker{\boxb}$ to $\boxb u =f$ if and only if $f\bot \Ker{\boxb}$. Furthermore if $f\in S^k$ then $u \in S^{k;2}$ and there is a uniform estimate
\[ \|u \|_{k;2} \leqc \|f \|_k. \]
Alternatively phrased, $\rng{\boxb} = \Ker{\boxb}^\bot$ and for all $k$,
\[ \boxb \colon S^{k;2} \cap \Ker{\boxb}^\bot \cap \dom{\boxb} \longrightarrow \Ker{\boxb}^\bot \cap S^k \]
is an isomorphism.
\end{thms}

The precise definition of the spaces and norms used here is given in  \rfS{MN}. We mention here that this theorem encodes exact regularity of solutions in the Folland-Stein spaces in all directions. Furthermore we obtain a full gain of two Folland-Stein derivatives for all directions in the interior and in directions tangent to the foliation at the boundary. In particular, hypoellipticity at the boundary for the canonical solution to $\boxb$ is implied.

An important application of the \db-Neumann problem is to solving the inhomogeneous tangential Cauchy-Riemann equation. Fortunately, in this instance we can partially remove condition (A2).  Our main theorem yields an existence and regularity theory for this problem. Our result is as follows:

\begin{thms} \label{B}
Let $\Omega$ be a smoothly bounded domain in a strictly pseudoconvex pseudohermitian manifold $(M,\eta)$ of dimension $2n+1$ with $n \geq 3$ such that $\Omega$ has noncharacteristic boundary and $(\Omega,w,\eta)$ satisfies  (A1). Then when the Kohn Laplacian is known to have closed range, the system 
\[  \db \varphi = \vs, \qquad \db \vs =0 \]
is solvable for $\varphi \in L^2(\Omega)$ if and and only if $\vs \perp \Ker{\boxb}$. Furthermore there is a closed complement $L$ to $\Ker{\db}$ in $L^2(\Omega)$ such that if $\vs \in S^k$ then $\varphi$ can be uniquely chosen to lie within $S^{k;1} \cap L$ and there is a uniform estimate 
\[ \|\varphi\|_{k;1} \leqc \|\vs\|_k.\]
Equivalently, for $k \geq 0$
\[ \db\colon L \cap S^{k;1} \longrightarrow \Ker{\db} \cap \Ker{\boxb}^\bot\]
is an isomorphism.

If (A2) also holds, then we may may choose $L = \Ker{\db}^\bot$.
\end{thms}

Again the precise definitions of the spaces involved is put off until \rfS{MN}. However we note that this encodes exact regularity in the weighted Folland-Stein spaces with a slight  gain in directions tangential to the foliation. This is sufficient to establish solutions globally smooth up the boundary when $\vs$ is itself smooth.

Considering the conditions of Theorem \ref{A}, it becomes important to understand when the Kohn Laplacian has closed range as an unbounded operator on $L^2$. As a partial answer to this question, we shall prove the following theorem:

\begin{thms}\label{C}
If $1 \leq q \leq n-2$, $(\Omega,\eta,w)$ satisfies  (A1), $\Omega$ has noncharacteristic boundary and all boundary leaves of the foliation by $w$ have zero Kohn-Rossi cohomology in degree $(0,q)$ then $\boxb$ is a Fredholm operator on $L^2(\Omega)$. Furthermore if (A2) holds, then $\boxb$ is hypoelliptic up to the boundary.
\end{thms}

We can combine this with the following known result

\bgT{noco}
If $M$ is a compact strictly pseudoconvex manifold of dimension $2n-1$ $n\geq 3$ such that $M$ is embedded in $\cn{n}$ or $\cn{n+1}$, then $M$ has no Kohn-Rossi cohomology is degrees $(0,q)$ for $1 \leq q \leq n-2$. 
\enT

For the \cn{n} case, pseudoconvexity is sufficient and the result is easily derived from Theorem 9.4.2 in \cite{Shaw:B}. For the other case, work by Harvey-Lawson \cite{Har-Law} shows that $M$ bounds a variety in $\cn{n+1}$ with isolated singular points. By dimension count these singularities are of hypersurface type. Yau \cite{Yau} then computed the Kohn-Rossi cohomology explicitly in terms of the moduli spaces of the singularities, in particular showing that it vanishes in degrees $(0,q)$ with $1 \leq q \leq n-2$.

As a corollary of these results, we establish the following theorem:

\begin{thms}\label{D}
If $M$ is a strictly pseudoconvex hypersurface in $\cn{n+1}$ with $n \geq 3$ and $\Omega$ a smooth compact domain with noncharacteristic boundary satisfying (A1) then $\boxb$ is a Fredholm  operator on $L^2_{(p,q)}(\Omega)$ for $1 \leq q \leq n-2$. 
\end{thms}

Furthermore, Theorem \ref{B} is of greatest practical use in circumstances where $\boxb$ is not only Fredholm but actually injective. For then we have a very simple criteria to check for solvability of $\db$.  We combine our results with earlier work by Shaw \cite{Shaw:P} to establish the following:

\begin{thms}\label{E}
Let $M$ be a strictly pseudoconvex pseudohermitian manifold embedded as a hypersurface in \cn{n+1} with defining function $r$. Let  $ \Omega =M \cap \{ \varrho < 0 \}$ be a bounded domain in $M$ with smooth, strictly convex defining function $\varrho = \varrho(z^1, \bt{z}^1)$.  Suppose also that  $1\leq q \leq n-2$ and $dr \wedge dz^1 \wedge d\bt{z}^1 \ne 0$ on $\partial \Omega$ ( i.e. that $\Omega$ has non-characteristic boundary). Then for any $(p,q)$-form $\varphi \in \Cc{\Omega}$ such that $\db \varphi =0$  there exists $(p,q-1)$-form $u \in \Cc{\Omega}$ such that $\db u = \varphi$.
\end{thms}

\section{Basic Definitions}\setS{PS}
 
 A pseudohermitian manifold consists of a triple $(M,\eta,J)$ where $M$ is a smooth, $2n+1$ dimensional real manifold, $\eta$ is a non-vanishing $1$-form on $M$ and $J:H \to H$ is a smooth bundle map on $H:=\ker{\eta}$ with $J^2=-1$. It is further assumed that the integrablity condition that $\crT$ and $\acrT$, the $+i$ and $-i$ eigenspaces of $J$ in the complexification of $H$ respectively, are involutive. Thus a pseudohermitian manifold can be considered as a codimension $1$ CR manifold together with a fixed, global contact form. 
 
 The Levi form for $M$ is the bilinear form $(X,Y) \mapsto d\eta(X,JY)$ on $H$. The structure is said to be strictly pseudoconvex if the Levi form is positive definite everywhere. In this case, there is a unique global vector field $T$ known as the characteristic field satisfying $T \lrcorner \eta =1$ and $T \lrcorner d\eta =0$. Thus for strictly pseudoconvex pseudohermitian manifolds, we can naturally extend $J$ by setting $JT=0$ and create a canonical metric
 \[ h_\eta(X,Y) = d\eta(X,J\Bt{Y})+\eta(X)\eta(\Bt{Y})\]
 which is Riemannian on $TM$ and Hermitian on $\cn{}TM$. The complexified tangent bundle then orthogonal decomposes as 
 \[ \cn{}TM = \crT \oplus \acrT \oplus \cn{}T.\]
On a strictly pseudoconvex pseudohermitian manifold there is a canonical connection. This allows us to intrinsically define a variety of Sobolev type spaces in addition to providing a useful computational tool.

\bgL{Connection}
If $(M,J,\eta)$ is strictly pseudoconvex, there is a unique connection $\nabla$ on $(M,J,\eta)$ that is compatible with the pseudohermitian structure in the sense that $H$, $T$, $J$ and $d\eta$ are all parallel and the torsion satisfies
\begin{align*}
\text{Tor}(X,Y)&=d\eta(X,Y)T\\
\text{Tor}(T,JX)&= -J\text{Tor}(T,X)
\end{align*}
\enL

This formulation of the connection was developed by Tanaka \cite{Tanaka}. An alternative formulation in terms of a coframe was independently derived by Webster  \cite{Webster}.  Tanaka's proof is constructive and is based upon the useful formula
\bgE{Connection}
\nabla_{X\dupp} Y\upp = [X\dupp,Y\upp]\upp.
\enE
where $X\upp$, $X\dupp$ denote the orthogonal projection of a vector field $X$ to $\crT$, $\acrT$ respectively. As the connection and metric naturally extend to other tensor bundles over $M$, we can instrinsically define $L^2$-Sobolev spaces via the norms
 \bgE{Sobolev} \|\varphi\|_{H^j} = \sum\limits_{k \leq j}\| \nabla^k \varphi\|_{L_\eta^2}. \enE
Unfortunately, these spaces do not provide optimal results for the analysis of the Kohn Laplacian as it is not fully elliptic; it's only first order in the characteristic direction. The Folland-Stein spaces were introduced in \cite{Folland:H} to provide more refined regularity results. They are defined intrinsically. For a differential form $\varphi$, we can decompose \[\nabla \varphi = \nabla \upp \varphi + \nabla \dupp \varphi + \nabla_T  \varphi \otimes \eta.\]   Here we define  $\nabla\upp$ by $\nabla\upp \varphi (\cdot, X)=\nabla \varphi (\cdot, X\upp)$ and $\nabla\dupp$ by $\nabla\dupp \varphi(\cdot, X)=\nabla \varphi(\cdot,X\dupp)$.
Now set 
\[ \nabla_H  \varphi := \nabla \upp \varphi + \nabla\dupp \varphi= \nabla \varphi - \nabla_{T} \varphi \otimes \eta.\]
The Folland-Stein spaces $S^j$ are now defined in from the norms
with norms
 \bgE{FollandStein} \| \varphi\|_j = \sum\limits_{k \leq j} \left\|\left( \nabla_H \right)^k \varphi\right\|_{L_\eta^2}. \enE

For the remainder of this section we suppose $(M,J,\eta)$ is a strictly pseudoconvex pseudohermitian manifold. Set $\Lambda_\eta^{0,1}M=\{ \varphi \in \cn{}T^*M: \varphi=0 \text{ on } \crT \oplus \cn{}T\}$ and let $\Lambda_\eta^{1,0}M$ be the orthogonal complement to $\Lambda_\eta^{0,1}M$ in $\cn{}T^*M$. It should be stressed that these definitions are asymmetric. We extend to higher degree forms by setting $\Lambda_\eta^{p,q}M = \Lambda^p\left( \Lambda_\eta^{1,0} M \right)\otimes \Lambda^q \left(\Lambda_\eta^{0,1}M\right)$. The space of degree $k$ complex covector fields on $M$ then admits the following orthogonal decomposition $\cn{}\Lambda^{k} M = \bigoplus\limits_{p+q=k} \Lambda_\eta^{p,q}M$. Denote the orthogonal projection $\cn{}\Lambda^{p+q}M \to \Lambda_\eta^{p,q}M$ by $\pi_\eta^{p,q}$ and define
\bgE{concrete}
\db = \pi_\eta^{p,q+1} \circ d.
\enE
Then $\db$ maps $C^\infty(\Lambda_\eta^{p,q}) \to C^\infty(\Lambda_\eta^{p,q+1})$. It should be remarked that this definition depends upon the pseudohermitian structure and is not canonical for the underlying CR structure. Using the language of holomorphic vector bundles and quotients, it is possible to construct an operator depending solely on the CR structure that reduces to our definition once a pseudohermitian form is chosen. However, for the purposes of this paper the concrete version offered here will suffice. It is also easy to see that for a smooth $(0,q)$-form $\varphi$,
\bgE{Computedb}
 \db \varphi = (-1)^{q}(q+1)\textup{Alt } (\nabla\dupp \varphi).
\enE

It follows immediately from the definitions that $\db \circ \db =0$. Thus $\db$ defines a complex of differential forms on $M$. The associated cohomology is known as the Kohn-Rossi cohomology and is denoted by $\mathcal{H}^{p,q}(M)$. A key tool for studying these groups is the Kohn Laplacian.
\bgD{Kohn}
If $(M,J,\eta)$ is a strictly pseudoconvex structure then the formal Kohn Laplacian is defined by
\[ \boxb = \db \vtb + \vtb \db \]
where \vtb denotes the formal adjoint of \db with respect to the canonical $L^2$ inner product on $(M,J,\eta)$.
\enD
On compact manifolds this operator is well-understood. For example we have the following theorem due to Folland  and Stein \cite{Folland:H}.
\bgT{compact} Let $(M,J,\eta)$ be a compact, strictly pseudoconvex pseudohermitian manifold of dimension $2n-1$. 
\bgEn{\alph}
\etem If $1 \leq q \leq n-2$ then $\boxb$ is a self-adjoint, Fredholm operator on $L^2(\Lambda_\eta^{p,q}M)$ and there is an orthogonal decomposition
\[ \begin{split}L^2(\Lambda_\eta^{p,q}M) &= \rng{\boxb} \oplus \Ker{\boxb}\\
&= \rng{\db} \oplus \rng{\vtb} \oplus \Ker{\boxb}.\end{split}\]
The operator $\boxb$ is subelliptic. Therefore $\Ker{\boxb}$ is finite dimensional  and $(1+\boxb)^{-1}$ is a compact, bounded operator on $L^2$. The cohomology group $\mathcal{H}^{p,q}(M) \cong \Ker{\boxb}$. Furthermore the operator $1+\boxb$ is an isomorphism from $S^{k+2}(\Lambda^{p,q}_\eta M)$ to $S^k(\Lambda^{p,q}_\eta M)$ for all $k \geq 0$. 
\etem If $q=0$ then $\boxb$ is self-adjoint and has closed range as an operator on $L^2(\Lambda_\eta^{p,0}M)$ and there is an orthogonal decomposition
\[ L^2(\Lambda_\eta^{p,0}M) = \rng{\boxb} \oplus \Ker{\boxb}.\]
\enEn
\enT

\bgR{dimension}
The reader is cautioned that the dimension assumed in the theorem is $2n-1$ rather than the $2n+1$ used earlier. When we are adding the supposition that the manifold is compact we shall always adopt this drop of dimension, whereas if we are not presupposing compactness we shall continue to use $2n+1$. This is to ensure compatibility of results when we are working with a foliation of a domain by compact CR manifolds of codimension $2$.
\enR

\bgR{OtherP}
Throughout the literature, most computations and arguments concerning $\db$ are conducted under the assumption that $p=0$. When the bundle $\Lambda^{1,0}_\eta M$ is holomorphically trivial, it is easy to pass to the general case. Since $\cn{}TM/\acrT $ is a holomorphic bundle, local results always extend simply to the case $p>0$.  If $M^{2n-1}$  is globally CR embeddable as a hypersurface in $\cn{n}$ then the embedding functions  provide a global trivialisation. \enR

For a compact manifold, the formal and actual $L^2$ adjoints of \db are equal. However for manifolds with boundaries the issue of boundary conditions arises and we must make a subtly different definition to recover self-adjointness for the operator.

Suppose $\Omega$ is a bounded open set in $M$ with smooth boundary. We can restrict our complexes to forms defined on $\Omega$. We extend $\db$ to its maximal $L^2$  closure, also denoted \db, and define $\dbs$ to be the $L^2$-adjoint of this extended operator. We then define the Kohn Laplacian for $\Omega$ by
\[\dom{\boxb} = \{\varphi \in \dom{\db} \cap \dom{\dbs}: \text{ $\db \varphi \in \dom{\dbs}$ and $\dbs \varphi \in \dom{\db}$} \}\]
and for  $\varphi \in \dom{\boxb}$
\[ \boxb \varphi = \db  \dbs \varphi + \dbs  \db \varphi.\] 
This operator is  self-adjoint as an operator on $L^2(\Lambda_\eta^{p,q}\Omega)$ (see \cite{Shaw:B}). For forms contained in $\dom{\boxb}$ the operator agrees with the formal version defined above. The $\db$-Neumann problem on $\Omega$  is then to decide when the equation $\boxb u =f$ on $\Omega$ can be solved for $u \in \dom{\boxb}$  and obtain optimal regularity results. It is worth emphasizing that there are boundary constraints on any solution $u$ for it to lie in $\dom{\boxb}$.  From the view point of CR geometry as opposed to pseudohermitian geometry we would also be free to choose an appropriate $\eta$. 

The analysis of this problem is difficult for several reasons. The operator is not elliptic as it has only limited control over the characteristic direction. However, the Folland-Stein spaces were constructed to address precisely this.  The characteristic vector field $T$ can be written as a commutation of vector fields from $H$. Thus $T$ is second order as an operator in the Folland-Stein setting. Although \boxb is only subelliptic (see \cite{Shaw:B}), it is fully elliptic in the Folland-Stein directions. 
A second problem  is that the boundary conditions are non-coercive in a sense to be made more precise later. There is also a third problem related to the geometry of the boundary \dO.

\bgD{CharPoint}
A point $x \in \dO$ is a characteristic point for $\Omega$ if the boundary is tangent to the distribution $H$ at $x$, i.e. $T_x\dO =H_x$.
\enD

At all non-characteristic points the tangent space to $\dO$ intersects $H$ transversely with codimension $1$. Thus at characteristic points there is a jump in the dimension of this intersection. This phenomenon makes obtaining $L^2$ estimates difficult near these points. 

All positive results for this problem have required strict conditions on the geometry of the boundary of $\Omega$ and have returned non-sharp boundary regularity. See for example the work by Ricardo Diaz \cite{Diaz} or Mei-Chi Shaw \cite{Shaw:P}. In particular very little is known about regularity when the domain possesses characteristic points.

At times throughout this paper we shall impose various conditions on the domain $\Omega$, the pseudohermitian form $\eta$ and a CR function $w$. In particular, we shall always be supposing (A1).

\begin{description}
\item[(A1)] $\Omega$ is a precompact open set open set with smooth defining function $\varrho$ depending solely on the real and imaginary parts of $w$.
\item[(A2)] Near the boundary $\partial \Omega$, the function $H_\bt{w} w$ is constant on the level sets of $w$.
\item[(A3)] The function $H_\bt{w} w$ is globally constant on the level sets of $w$.
\item[(B)]   $dw \wedge d\bt{w} \ne 0$ on $\Bt{\Omega}$.
\end{description}
The pseudohamiltonian vector field $H_\bt{w}$ is defined in \rfS{NL}. However, we note here that whenever $\Omega$ satisfies (A1) and is noncharacteristic we can choose a pseudohermitian form such that (A2) holds. Condition (B) is a more stringent requirement than noncharacteristic boundary, forcing the level sets of $w$ to be nondegenerate globally on $\Omega$.  If (B) holds then the pseudohermitian form can be chosen so that (A3) holds.

\section{Normalization}\setS{NL}

Frequently, it is the underlying CR manifold that is of interest rather than the specific pseudohermitian structure. One advantage of this is that we can scale the pseudohermitian form to simplify computations. Throughout this section we assume that $(M,J,\eta)$ is a $2n+1$-dimensional strictly pseudoconvex pseudohermitian manifold and that $w$ is a CR function $M$.

\bgD{Hf}
The pseudohamiltonian field for a smooth function $x$, $H_x$,is defined by
\[ \eta(H_x) =0, \qquad H_x \lrcorner d\eta = d_b x \quad (:= dx -(Tx)\eta) \]
\enD
We note in passing that
\[H_x v = d\eta(H_v,H_x) = -d\eta(H_x,H_v) = -H_v x.\]
Pseudohamiltonian  fields are the key ingredient to understanding how the non-horizontal vector field bracket structure changes as we rescale the pseudohermitian form.
\bgL{Tu}
The characteristic field for the pseudohermitian form $\eta^\au = e^x \eta$ is
\[ T^\au = e^{-x} (T + H_x) \]
\enL

\pf This is just a matter of computation
\[ T\lrcorner d\eta^\au = T_\lrcorner e^x du \wedge \eta + e^x T \lrcorner d\eta = -e^x d_b x,\]
\[ H_x \lrcorner d\eta^\au = e^x (H_x u) \eta + e^x H_x \lrcorner d\eta = e^x d_b x. \]
\epf

One complication when comparing operators on rescaled pseudohermitian structures is in the differing presentations of $(p,q)$-forms. Since the spaces $\Lambda_\eta^{p,q}$ depended on orthogonal projections for the metric induced by $\eta$, we get get a different space for $\eta^\au = e^x \eta$ whenever the smooth real-valued function $x$ is not identically zero. However we can introduce operators $\mu_x \colon \Lambda_\eta^{0,1} \to \Lambda_{\eta^\au}^{0,1}$ by
\bgE{mu} \mu_x \varphi = \varphi - \varphi(T^\au) \eta^\au. \enE
We can immediately extend this to $(0,q)$-forms by declaring $\mu_x ( \varphi \wedge \psi ) = \mu_x \varphi \wedge \mu_x \psi$. 

\bgL{Hu}
\begin{align*}
d^\au_b f &= d_b f -(H_x f) \eta \\
H^\au_f  &= e^{-x} H_f\\ 
\end{align*}
\enL

\pf Again, we compute
\begin{align*}
d^\au_b f &= df -(T^\au f) \eta^\au = df - (Tf)\eta -(H_x f) \eta = d_b f - (H_x f) \eta,\\
 H_f \lrcorner d\eta^\au &= e^x (H_f x) \eta + e^x d_b f = e^x \left( df-(Tf)\eta -(H_x f) \eta \right) =e^x d_b^\au f.\\
\end{align*}
\epf

The heart of method presented in this paper is understanding how the Kohn Laplacian behaves with respect to a foliation by level sets of the CR function, $w$. A key step shall be decomposing the operator into pieces tangent and transverse to the foliation. Accordingly, we shall now describe the canonical vector fields from which the transverse pieces will be constructed.
 
\bgD{Y}
\[ Y= H_\bt{w}, \qquad Y^\au = H^\au_\bt{w} \] 
\enD

\bgL{Ynorm}
\[ \aip{Y^\au}{Y^\au}{\au} = -i\Bt{Y}^\au \bt{w} = -ie^{-x}\Bt{Y} \bt{w} \]
\enL

\pf First compute
\[ \aip{Y^\au}{Y^\au}{\au} = d\eta^\au(Y^\au,J\Bt{Y}^\au)= -i (\db[\au] \bt{w})(\Bt{Y}^\au) = -i \Bt{Y}^\au \bt{w}= -ie^{-x} \Bt{Y} \bt{w} \]
\epf

From this we see that  away from the characteristic locus of $w$
\[ \chs =\left\{ p\in M \colon  \db w(p) =0\right\} =\left\{ p \in M \colon (dw\wedge d\bt{w}) _{|p} =0 \right\}  \]
we can then fix a canonical pseudohermitian form for $w$ by fixing $\Theta = \eta^\au$ to be the form such that
\[ \aip{Y^\au}{Y^\au}{\au} = 1 \quad \text{ or equivalently }\quad  \Bt{Y}^\au \bt{w} = i.\]
It should be noted that for a smooth domain $\Omega$ satisfying (A1), the characteristic set of $\dO$ is exactly $\dO \cap \chs$. Therefore on a noncharacteristic domain we can normalize the pseudohermitian form so that (A2) holds. If (B) holds, then $\chs \cap \overline{\Omega} = \emptyset$ and the pseudohermitian form can be normalized so that (A3) holds.


\section{Decomposition of the Kohn Laplacian}\setS{CP2}

In this section we shall construct the decomposition of the Kohn Laplacian on which we shall base our regularity results. We suppose $w$ is a fixed CR function on the strictly pseudoconvex pseudohermitian manifold $(M^{2n+1},\eta,J)$ and $\Omega$ is a smoothly bounded precompact open set such that $(\Omega, \eta ,w)$ satisfy (A1) and (A3).  Thus the pseudohermitian form $\eta$ has the property that 
\[ \left| \Bt{Y} \right|^2 =-i \Bt{Y} \bt{w} = e^{2\nu} >0 \]
where $\nu = \nu(w,\bt{w})$. 

We shall adopt the convention that latin indices run from $0$ to $n-1$ and greek indices run from $1$ to $n-1$. To fit in with this convention, we set \[Z_0 = e^{-\nu} Y.\] We shall work with a local orthonormal frame $Z_0, Z_1, \dots , Z_{n-1}$ for $\crT$ with the property that each $Z_\ua$ is tangent to the level sets of $w$ and $Z_0$ is defined as above. The dual frame will be denoted $\theta^0, \theta^1 , \dots ,\theta^{n-1}$.  Note that this implies that $Z_\ua \nu = 0 = Z_\ba \nu$.

If $\wh{p}$ is the leaf of the foliation by $w$ containing $p$, we let $\iota\colon \wh{p} \hookrightarrow M$ be the inclusion map and consider the pseudohermitian form $\wh{\eta} = \iota^* \eta$ with $\wh{J}= J_{|T\wh{p}}$.

\bgD{Connection} Let $\om^k_m$ and $\wom^\ua_\ub$ are the connection 1-forms for the Tanaka-Webster connections for $(M,\eta,J)$, $(\wh{p},\wh{\eta},\wh{J})$ respectively associated to the frames $\{ Z_j \}$ and $\{Z_\ua\}$. Thus
\[ \nabla Z_k = \om_k^m \otimes Z_m, \qquad \wh{\nabla} Z_\ua = \wom_\ua^\ub \otimes Z_\ub\]
We also use $\G^m_{jk} =\om^m_k(Z_j)$ and $\wG^\ua_{\ub \uc} = \wom^\ua_{\uc}(Z_\ub)$.
\enD

Important computational tools are the following structural equations for Tanaka-Webster connection \cite{Webster}.

\bgE{TW}
\begin{split}
d\eta &= ih_{j \bt{k}} \theta^j \wedge \theta^\bt{k} = i \theta^\ua \wedge \theta^\ba + i \theta^0 \wedge \theta^\bt{0} \\
dh_{j \bt{k} }&= h_{m \bt{k} } \om^m_j + h_{j \bt{m}} \om^\bt{m}_{\bt{k}}  = \om^k_j + \om^\bt{j}_\bt{k}\\
d\theta^\bt{k} &= \theta^\bt{m} \wedge \om_\bt{m}^\bt{k} + \eta \wedge A^\bt{k}_m \theta^m\\
\end{split}
\enE
where $h_{j \bt{k}}$ are the components of the Levi metric. Our chosen frame is orthonormal, so  $h_{j \bt{k} } = \begin{cases} 1, & i=j \\ 0, & i \ne j  \end{cases}$.

Using these we can now start to explore the relationships between the connection on $M$ and those on the foliating leaves $\wh{p}$.

\bgL{tangentT}
The characteristic vector field for $(\wh{p},\wh{\eta}, \wh{J})$ is 
\[ \wh{T} = T +aY +\bt{a} \Bt{Y} \] 
where $a= -ie^{-2\nu} T w$ 
\enL

\pf
Clearly $\wh{T} \lrcorner \wh{\eta} =1$. Now since
\[ d\eta = i \theta^\ua \wedge \theta^\ba + i \theta^0 \wedge  \theta^\bt{0}\]
we have 
\[ d\wh{\eta} = \iota^* d\eta = i \theta^\ua \wedge \theta^\ba.\]
Thus 
\[ \wh{T} \lrcorner d\wh{\eta} = 0.\]
The final requirement is that $\wh{T}$ is real and is tangent to $\wh{p}$. But $\wh{T}$ is clearly real and the tangency condition follows from 
\[ \wh{T} w = T w +aYw = T w  -ia = 0\]
and 
\[ \wh{T} \bt{w} = \Bt{\wh{T}  w} =0.\]
\epf

\bgC{w0}
\[ \iota^* \theta^0 = a \wh{\eta}, \qquad \iota^* \theta^\bt{0}= \bt{a} \wh{\eta}\]
\enC

\bgL{WTM}
The Christoffel symbols $\G^m_{jk}$ satisfy the following properties
\begin{align}
\G^k_{\bullet j } &= -\G^\bt{j}_{\bullet \bt{k}   }  \\
\G^j_{\bullet \bt{k}} &=0\\  
\G^0_{ \bb \ua} &= ia \delta_{\ua \bb}\\
 \G^0_{ j k} &=\G^0_{ k j}\\
\G^\bt{0}_{ \ua \bt{0}} &=  \G^\ba_{ 0 \bt{0} }=\G^{\ua}_{00}=\G^\bt{0}_{0 \ba }=0
 \end{align}
\enL

\pf
The first is just a defining property of the Tanaka-Webster connection
\[ dh_{j \bt{k}} = h_{j \bt{m}} \om^\bt{m}_\bt{k} + h_{m \bt{k}}\om^m_j.\]
The second is just that $\acrT$ is parallel.

For the third we note
\[ \G^0_{\bb \ua } = \G^0_{ \bb \ua} - \G^0_{\ua \bb } = \theta^0 ([Z_\bb,Z_\ua]) \]
Now $[Z_\ua,Z_\bb]$ must be tangent to $\wh{p}$  and $\eta([Z_\ua,Z_\bb]) = -i \delta_{\ua \bb}$. \rfL{tangentT} then implies that $\theta^0([Z_\bb,Z_\ua])= ia \delta_{\ua \bb}$. 

For the fourth note that an easy consequence of $|Y|=e^{\nu}$  depending on $w$ and $\bt{w}$ alone is that all Lie brackets of the form $[Z_j,Z_\ua]$, $[Z_0, Z_\ba]$ and their conjugates, while horizontal,  have no $Z_0, Z_\bt{0}$ components. Therefore
\begin{align*}
 \G^0_{ jk} &= d\theta^0 (Z_j,Z_k) = \theta^\bt{0} ([Z_j,Z_k]) =0
 \end{align*}
 For the fifth, again we note
  \begin{align*}
 \G^{\ba}_{0 \bt{0} } &= -\G^0_{0 \ua } = -\G^0_{ \ua 0}  = \G^\bt{0}_{\ua \bt{0} }  = \theta^\bt{0} ([Z_\ua, Z_\bt{0}]) =0  \\
  \G^\ua_{0 \bt{0} } &= -\G^\bt{0}_{0 \ba }  = -\theta^\bt{0}([Z_0,Z_\ba]) =0 
  \end{align*}
\epf

\bgP{WChristoffel}
The connection and torsion forms for the connections $\nabla$ and $\wh{\nabla}$ are related as follows:
\begin{align*}
\wG^\ua_{ \uc \ub} &= \G^\ua_{\uc \ub }\\
\wG^\ua_{ \bc \ub} &= \G^\ua_{ \bc \ub}\\
\wG^\ua_{T \ub } &= \G^\ua_{T \ub } +a \G^\ua_{0 \ub } +\bt{a}\G^\ua_{\bt{0} \ub } - a\G^\ua_{\ub 0} \\
&= \G^\ua_{T \ub } +a \G^\ua_{0 \ub } +\bt{a}\G^\ua_{ \bt{0} \ub}  -i|a|^2 \delta_{\ua \ub}\\
\wh{A}^\ua_\bb &= A^\ua_\bb +a \G^\ua_{ \bb 0}.
\end{align*}
\enP

\pf   All that is required is to check that these forms obey the structural equations. This routine computation is left to the reader.

\epf

When acting on $(0,q)$-forms, it follows easily from the definitions that he tangential CR operator and its formal adjoint can be expressed in terms of the Tanaka-Webster connection as
\bgE{CR}
 \db \varphi = \theta^\bt{k} \wedge \nabla_\bt{k} \varphi, \qquad \vartheta_b \varphi =- \theta^k \vee \nabla_k \varphi.
\enE
Thus we can use our comparison of the connections $\nabla$ and $\wh{\nabla}$ to break these operators down in to pieces tangent and transverse to the folation.

\bgL{Nabla}
\begin{align}
(\nabla - \wh{\nabla} )_\bb \theta^\ba &= ia \delta_{\ba \bb} \theta^\bt{0} \\
(\nabla - \wh{\nabla} )_\ub \varphi &= \G^0_{ \ub \ua}    \theta^\bt{0} \wedge \theta^\ua \vee \varphi \\
\theta^0 \vee \nabla_0 \theta^\ba &= 0
\end{align}
\enL

\pf Again, we compute
\begin{align*}
(\nabla - \wh{\nabla} )_\bb \theta^\ba &= -\G^\ba_{ \bb \bt{k}} \theta^\bt{k} +\wG^\ba_{ \bb \bc} \theta^\bc \\
&= -\G^\ba_{ \bb \bt{0}} \theta^\bt{0} = \G^0_{ \bb \ua} \theta^\bt{0} \\
&= ia \delta_{\ba \bb} \theta^\bt{0}
\end{align*}

\begin{align*}
(\nabla - \wh{\nabla} )_\ub \theta^\ba &= -\G^\ba_{ \ub \bt{k}} \theta^\bt{k} +\wG^\ba_{ \ub \bc} \theta^\bc \\
&= -\G^\ba_{ \ub \bt{0}} \theta^\bt{0} = \G^0_{\ub \ua } \theta^\bt{0} \\
&= -\G^0_{ \ub \ua} \theta^\ua \vee \theta^\bt{0} \wedge \theta^\ba
\end{align*}

\begin{align*}
\theta^0 \vee \nabla_Y \theta^\ba &= \theta^0 \vee (- \G^\ba_{0 \bt{k} } \theta^\bt{k} )\\
&= - \G^\ba_{0\bt{0} } = 0
\end{align*}
\epf

\bgD{Dq}
For non-negative integers $q$, we define $\oD_q$ on smooth $(0,q)$-forms $\varphi$ by
\[ \oD_q  \varphi= \nabla_\bt{0}  \varphi- iaq \varphi \]
and extend as a maximal closed operator on $L^2$. 
\enD

\bgL{Dq*}
The formal adjoint, $D^{\#}_q$ of $\oD_q$ is given by 
\[ D^{\#}_q \varphi =   \G^\bt{0}_{0 \bt{0} } \varphi - i\bt{a}(n-q-1)\varphi  - \nabla_0 \varphi   \]
\enL

\pf
From a standard divergence theorem argument $\nabla_\bt{0}^*  = -\text{div }Y  - \nabla_0$. Now
\begin{align*}
 \text{div }Y &= \G^\bullet_{ \bullet 0} = \G^k_{ k 0} =\G^{\bt{0}}_{k \bt{k} } \\
&= -\G^\bt{0}_{0 \bt{0} } + i\bt{a}(n-1) 
\end{align*}
\epf

\bgD{V} We define the operators  $\hdb$, $\hdbs$ and $\hboxb$ by
\begin{align*} \hdb \varphi &= \theta^\ba \wedge \wh{\nabla}_\ba \varphi\\   \hdbs \varphi &=- \theta^\ua \vee \wh{\nabla} _\ua \varphi\\ \hboxb  \varphi &=  \hdb \hdbs  \varphi+ \hdbs \hdb \varphi  \end{align*}
\enD
Thus $\hdbs$ is the formal adjoint of $\db$ and $\hboxb$ acts as the Kohn Laplacian on each foliating leaf.

\bgL{db} Suppose $\varphi \in \Ci$ is a $(0,q)$-form such that $\oY \lrcorner \varphi =0$ then
\begin{align} \db \varphi &= \hdb \varphi +  \theta^\bt{0} \wedge \oD_q \varphi\\
\db \left( \theta^\bt{0} \wedge \varphi \right)&= -\theta^\bt{0} \wedge   \hdb \varphi \\
\fdbs \varphi &= \hdbs \varphi \\
\fdbs \left( \theta^\bt{0} \wedge \varphi\right) &= D^\#_q \varphi - \theta^\bt{0} \wedge  \hdb  \varphi
    \end{align}
\enL

\pf Once more the proof is by computation. Since $\Bt{Y} \lrcorner \varphi=0$, the $(0,q)$-form $\varphi$ can be written as a linear combination of wedge products of the forms $\theta^\ba$. Thus 
\begin{align*}
\db \varphi  &= \theta^\bt{k} \wedge \nabla_\bt{k} ( \varphi )  = \theta^\ba \wedge \nabla_\ba (\varphi) + \theta^\bt{0} \wedge \nabla_\bt{0}(\varphi) \\
&= \hdb \varphi + \theta^\ba \wedge (\nabla -\wh{\nabla})_\ba (\varphi) + \theta^\bt{0} \wedge \nabla_\bt{0} \varphi\\
&= \hdb \varphi + \theta^\bt{0} \wedge (\nabla_\bt{0}-iaq) \varphi
\end{align*}
Now from the structural equations \rfE{TW} we see
\begin{align*}
d\theta^\bt{0} &= \theta^\bt{k} \wedge \om^\bt{0}_\bt{k} + \eta \wedge A^\bt{0}_k \theta^k\\
\end{align*}
and so
\[
 \db \theta^\bt{0}(Z_\ba,Z_\bb) =  \G^\bt{0}_{\bb \ba } -\G^\bt{0}_{\ba \bb} = 0 
\]
 whereas
 \[  \db \theta^\bt{0}(Z_\bt{0},Z_\ba) = \G^\bt{0}_{ \ba \bt{0} } -\G^\bt{0}_{ \bt{0}\ba} = \G^0_{ \ba 0} - \G^\ua_{  \bt{0} 0 } =  0.\]
 Thus
 \[ \db \theta^\bt{0} =  0  .\]
The second identity easily follows from this and the first. To prove the third, we compute similarly
 \begin{align*}
 \fdbs \varphi &= -\theta^k \vee \nabla_k \varphi =- \theta^0 \vee \nabla_0 \varphi - \theta^\ua \vee ( \nabla -\wh{\nabla})_\ua \varphi  + \hdbs \varphi \\
 &= 0  +\G^0_{\ub \ua} \theta^\ub \vee \theta^\ua \vee \theta^\bt{0} \wedge \varphi + \hdbs \varphi\\
 &= \hdbs \varphi.
 \end{align*}
 The fourth follows another similar computation
 \begin{align*}
 \fdbs \left( \theta^\bt{0} \wedge \varphi\right) &= -\theta^k \vee \nabla_k \left(\theta^\bt{0} \wedge\varphi\right)\\
 &= -\theta^k \vee \nabla_k ( \theta^\bt{0} ) \wedge \varphi -\theta^k \vee \theta^\bt{0} \wedge \nabla_k \varphi \\
 &=  \G^\bt{0}_{k \bt{m} } \theta^k \vee \theta^\bt{m} \wedge \varphi  -  \theta^0 \vee \theta^\bt{0} \wedge \nabla_0  \varphi  -\theta^\ua \vee \theta^\bt{0} \wedge \nabla_\ua \varphi \\
 &= \G^\bt{0}_{0\bt{0} } \varphi  + \theta^\bt{0} \wedge (\nu_\ua \theta^\ua) \vee \varphi -i\bt{a} \theta^\ua \vee \theta^\bt{a} \wedge  \varphi\\
 &\qquad  - (\nabla_0 \varphi)^\top + \theta^\bt{0} \wedge \theta^\ua \vee (\nabla - \wh{\nabla})_\ua \varphi - \theta^\bt{0} \wedge \hdbs \varphi\\
 &= \G^\bt{0}_{0\bt{0} } \varphi - i(n-1-q)\bt{a} \varphi - \nabla_0\varphi \\
 & \qquad - \theta^\bt{0} \wedge   - \hdbs \varphi .
\end{align*}
 
\epf

\bgC{decomposition}
For smooth $(0,q)$-forms $\varphi$ in the domain of $\boxb$ we have
\bgE{decompositionV} 
\begin{split}
\boxb \varphi &= \Big( \hboxb  + D^{\#}_q \oD_q \Big) \varphi^\top + \big[\hdb, D^{\#}_{q-1}\big]\varphi^\bot  \\
& \qquad + \theta^\bt{0} \wedge  \Big(   \hboxb +\oD_{q-1} D^{\#}_{q-1} \Big )\varphi^\bot + \theta^\bt{0} \wedge \big[ D_q, \hdbs\big] \varphi^\top
\end{split}
\enE 
\enC

This decomposition is dependent on  condition (A3), which greatly simplified the computations above. This condition will also aid commutation arguments in later sections. To express this concisely, we shall introduce some further terminology to be used throughout the paper.
\bgD{Notation}\hfill

\begin{itemize}
\item $\Psi\n[k]$ shall denote the set of all operators formed as the smooth compositions of $m$ with $0 \leq m \leq k$ covariant derivatives taken with respect to horizontal vector fields tangent to the foliating leaves. If $m=0$ the operator should be viewed as multiplication by a smooth function. 
\item $\Psi\n[{j,k}]$ refines this to the composition of at most $j$ $(1,0)$ vector fields and at most $k$ $(0,1)$ vector fields.
\item $\psi\n[k]$ (respectively $\psi\n[{j,k}]$) will indicate an element of $\Psi\n[k]$ (respectively $\Psi\n[{j,k}]$).
\item We shall use $S\n[k]$ and $S\n[{j,k}]$ in a fashion analogous to $\psi\n[k]$ and $\psi\n[j,k]$ but without the restriction that the vector fields be tangent to the foliating leaves.
\item $T\n[k]$ will likewise denote a generic $k$th order operator built out of $\wh{T}$ and zero order operators.
\item  If $q$ is understood, in this instance, we shall simply  write $\oD$ and $D$ in place of $\oD_q$ and $\oD^\#_q$.
\end{itemize}
\enD
For example we might express $[Z_\ua,Z_\bb]$ as any of  $\psi\n[2]$, $T\n +S\n$ or $T\n + \psi\n$.

\bgL{Properties}
Under Conditions (A1)  and (A3) the following properties hold

\bgEn{\alph}
\etem $[D, \hdb]$, $[\oD,\hdbs]$ are both in $\Psi\n$.
\item $[\psi_{0,1},\oD ] \in \Psi\n[{0,1}]$, $[\psi\n[{1,0}],D] \in \Psi\n[{1,0}]$
\etem $[\psi\n, D], [\psi\n,\oD] \in \Psi\n$.
\etem $[\hboxb, \oD], [\hboxb,D] \in \Psi\n[2]$
\enEn
\enL

\pf  We'll prove the first part of (a). The others are very similar.

Since $D = -\nabla_0 + \psi\n[0]$, it suffices to compute
\begin{align*}
[ \nabla_0 , \hdb ] &= [ \nabla_0 , \theta^\ba \wedge \nabla_\ba ] = \psi\n + \theta^\ba \wedge [\nabla_0 , \nabla_\ba] \\
&= \psi\n + e^{-\nu} \theta^\ba \wedge \nabla_{[ Y, Z_\ba ]} +e^{-\nu} \theta^\ba \wedge R(Y,Z_\ba) \\
&= \psi\n .
\end{align*}
The last equality follows as $[Y,Z_\ba]$ is a horizontal vector field which annihilates $w$ and $\bt{w}$ by (A3). Here $R$ represents the curvature endomorphism associated to $\nabla$.

\epf

\bgR{Operators}
Formally, we could now define operators
\[ \Pb =  \hboxb +\oD D, \qquad \Pt = \hboxb + D \oD \]
so that
\[ \boxb \varphi  = \Pt \varphi^\top + \psi\n[1] \varphi^\bot + \theta^\bt{0} \wedge \left( \Pb \varphi^\bot +\psi\n[1] \varphi^\top \right).\]
The plan is to study the properties of the subelliptic operators $\Pt$ and $\Pb$, then to absorb the error terms to obtain results for \boxb itself. To do this however, we shall in subsequent sections more carefully define these operators to take into account boundary conditions.
\enR

\section{Analysis on the foliation}\setS{FN}

In this section, we shall outline the properties of the Kohn Laplacians associated to the underlying foliation by compact CR manifolds . Throughout this section we shall assume that $\ud>0$ and $(\Omega, w,\eta)$ satisfy (A1), (A3) and (B). 

For convenience of notation we shall use $\leqc$ to indicate an inequality that holds up to multiplication by a positive constant that is independent of any function choices. For example $\|u \| \leqc \|Pu\|$ would indicate that there is a positive constant $C$ such that $\|u\| \leq C \|Pu\|$ for all $u$. If we wish to emphasize that the constant may depend on a parameter such as $\ud$, we will use the notation $\leqc[\ud]$ instead. 

Since we shall be distinguishing between transverse and tangential directions in our subsequent analysis, we define for $(0,q)$-forms the space of tangential smooth forms
\[ \Cit[\Omega] = \{ u \in \Ci[\Omega] \colon  \Bt{Y} \lrcorner u = 0 \} \]
with similar definitions for $L^{2,\top}(\Omega)$ and other functions spaces. We also note in passing that one consequence of condition (A3) is that $\nabla_\bt{0}$ preserves $\Cit[\Omega]$.

\bgD{L}
For $\ud>0$, we define the operator $\wh{P}_\ud$ on $\Cct{\Omega}$ by
\[ \wh{P}_\ud u = \begin{cases} \ud+ \hboxb, \quad &q>0\\ \ud+ \hboxb + \overline{\hboxb}, \quad &q=0  \end{cases}\]
where $\overline{\hboxb} u = \overline{ \hboxb \overline{u} }$.
\enD

Now for $\ud>0$, $\wh{P}_\ud$ is strictly positive and symmetric with respect to the $L^2$ inner product. 
Therefore we can create a new inner product
\[ \ip{u}{v}{P} := \ip{ \wh{P}_\ud u}{ v}{} \] and define $ \wh{S}^1$ to be the closure of $\Cct{\Omega}$ under the norm $u \mapsto  \ip{u}{u}{P}^{1/2}$. 

We set $\dom{\Lambda_\ud^2 u}$ to be the set
\[ \dom{\Lambda_\ud^2} = \left\{ u \in \wh{S}^1 \colon \left| \ip{u}{v}{P}\right| \leqc \| v\|_{L^2(\Omega)} \text{ for all $v \in \wh{S}^1$.} \right\}\]
and for $u \in \dom{\Lambda_\ud^2 u}$ we define $\Lambda^2 u$ by
\[ \ip{u}{v}{P} = \ip{\Lambda_\ud^2 u}{v}{},  \qquad v \in \wh{S}^1.\]
Then $\Lambda_\ud^2$ is a self-adjoint, strictly positive extension of $\wh{P}_\ud$. Thus we can define $\Lambda_\ud = ( \Lambda_\ud^2)^{1/2}$, a self-adjoint, strictly positive operator with $\dom{\Lambda_\ud} = \wh{S}^1$ and
\[ \ip{u}{v}{P} = \ip{\Lambda_\ud u}{\Lambda_\ud v}{}.\]

 \bgR{one}
 For the case $\ud=1$, we shall often simplify notation by setting $\Lambda = \Lambda_1$.
 \enR
 
Since $\Lambda_\ud$ is self-adjoint and strictly positive, we can define $\Lambda_\ud^s$ for all $s\in \rn{}$. 

\bgD{whS}
For $k>0$, we define $\wh{S}^k$ to be the closure of $\Cct{\Omega}$ in $L^{2,\top}(\Omega)$ under the inner product $\ip{u}{v}{k} := \ip{\Lambda^k u}{\Lambda^s k}{}$. 

For $k=-1$, we define $\wh{S}^{-1}$ to be the closure in $H^{-1,\top}(\Omega)$ of $\Cct{\Omega}$ under the inner product $\ip{u}{v}{-1}$.
\enD

Our first important result using these spaces is the following interpolation theorem

\bgL{Interpolation}
If $N \colon \Cct{\Omega} \to \Cct{\Omega}$ satisfies the following estimates
\[ \|Nf\|_{\wh{S}^3} \leq \ua \| f \|_{\wh{S}^1} , \qquad \|Nf\|_{\wh{S}^1} \leq \ub \| f\|_{\wh{S}^{-1}}\]
then it also satisfies the estimate
\[ \| Nf \|_{\wh{S}^2} \leq \max \{\ua,\ub\}  \| f\|_{L^{2,\top}(\Omega)}. \]
 \enL
 
 \pf
 The result essentially follows immediately from the classical interpolation theorem. The operator $N$ extends to a bounded operator $\wh{S}^1\to\wh{S}^3$ and $\wh{S}^{-1} \to \wh{S}^1$. Thus by interpolation (following the method and notation of Lion-Magenes \cite{LiMa})  $N$ is bounded from $[\wh{S}^{-1}, \wh{S}^1]_{1/2} \to  [\wh{S}^1,\wh{S}^3]_{1/2}.$ However, it is clear that $\Cct{\Omega}$ is densely contained in both these midpoint interpolation spaces and that on $\Cct{\Omega}$ the interpolation norms are equivalent to those of $L^{2,\top}(\Omega)$ and $\wh{S}^2$ respectively.
 
 \epf

The other result we shall need is the following regularity result.

\bgL{Regularity}
If $ \Lambda_\ud^2 u \in S^{k,\top} $ and $u \in S^{k-1,\top}$, then $u \in S^{k,\top}$.
\enL

\pf
On $\Cc{\Omega}$, we have that $\Lambda_\ud^2 = \wh{P}$. Thus $[ \Lambda_\ud^2 , \oY], [\Lambda_\ud, Y] = \psi\n[2]$. Now from standard regularity results on the foliating leaves, coupled with the smoothness of estimate coefficients in $w$, we see that if $\Lambda_\ud^2 u \in L^2$ then all $\psi\n[2] u \in L^{2,\top}(\Omega)$. Now formally by \rfL[CP2]{Properties}
\[ p_k(Y,\oY) \Lambda_\ud^2  =\Lambda_\ud^2 p_k(Y,\oY)+ \psi\n[2] p_{k-1} (Y,\oY).\]
Therefore an easy induction argument shows that $\Lambda_\ud^2 p_{k} (Y,\oY)u \in L^{2,\top}(\Omega)$. But this implies that $p_k(Y,\oY) u \in L^{2,\top}(\Omega)$. Derivatives tangent to the leaves are easily controlled using the regularity results of the leaves.

 \epf

\section{Transverse regularity and estimates}\setS{TV}

Throughout this section we shall assume that $0 \leq q \leq n-2$ and that $(\Omega,w,\eta)$ satisfy (A1), (A3) and (B). The set $\Cvt{\Omega}$ denotes the space of tangential smooth $(0,q)$ forms on $\Omega$ that vanish identically at all points of the boundary $\partial \Omega$. We also restrict the operators $\hdb$ and $\hdbs$ to tangential forms and extend to closed operators on $L^{2,\top}(\Omega)$.

\bgL{Density}
The space of compactly supported smooth forms, $\Coit[\Omega]$, is dense in 
\[\Bt{\mathcal{D}}= 
\dom{\hdb} \cap \dom{\hdbs} \cap \dom{\oD^*} \cap L^{2,\top}(\Omega).\]
in the graph norm. Note that if $q=0$ then we define $\dom{\hdbs}= L^{2,\top}(\Omega)$.
\enL

\pf
Since the derivative of $\hdb$ and $\hdbs$ are tangent to the boundary, the method of Lemma 4.3.2 (ii) in Chen-Shaw \cite{Shaw:B} can be applied.

\epf

\bgD{Q}
For $u \in \Cvt{\Omega}$ we define
\[\begin{split}
Q^\bot_{\e,\ud}(u,u) &= \ip{ \hdb u}{\hdb u}{} + \ip{\hdbs u}{\hdbs u}{} + \ip{Du}{Du}{}\\ & \qquad  + \ud\|u\|^2 +\e^2 \ip{ \wh{T} u}{\wh{T} u}{}\end{split} \]
and set $\Bt{\mathcal{D}}^\e$ to be the closure of $\Cvt{\Omega}$ under the $Q^\bot_{\e,1}$-norm. Note, the density lemma implies that $\Bt{\mathcal{D}}^0 = \Bt{\mathcal{D}}$.
\enD

\bgL{BasicEstimate}
For $\ud>0$, $\e \geq 0$ and $u \in  \Bt{\mathcal{D}}^\e$ there are uniform estimates
\begin{align*}
 \| u\|_1^2 +\e^2 \|u\|_{H^1}^2 &\leqc[\ud] Q_{\e,\ud}^\bot(u,u)\\
  \| u\|_1^2 +\e^2 \|u\|_{H^1}^2 &\leqc Q_{\e,\ud}^\bot(u,u) +\|u\| 
 \end{align*}
\enL

\pf  We'll prove $q=0$, the other case $q>0$ is easier as the estimate  \[\|\psi\n u\|^2 \leqc \| \hdb u\|^2 + \| \hdbs u\|^2 +\|u\|^2\] follows from standard estimates for the Kohn Laplacian on compact CR manifolds together with continuity of estimates in the parameter $w$.

Assume $q=0$. Note that in this case $\oD = e^{-\nu} \overline{Y}$. Now, essentially by definition $\| \psi\n[0,1] u \| \leqc \| \hdb u\|$. Next note that
\[ \left| \ip{\psi\n u}{u}{} \right| \leqc \left| \ip{ \psi\n[0,1] u}{u}{} \right| + \left| \ip{u}{\psi\n[0,1] u}{} \right| \leqc \| \hdb u\|^2 +\|u\|^2 \]
 Now
since $[Y,\oY] = -i \wh{T} + S\n$,
\begin{align*}
\| \oY^* u\|^2 &= \ip{ \oY \oY^* u}{u}{} = \ip{ \oY^* \oY u}{u}{} + \ip{[\oY,\oY^*]u}{u}{} \\
&= \| \oY u\|^2 +\ip{ i\wh{T} u}{u}{} + \ip{ (\psi\n[0] \oY  \psi\n[0]  \oY^* + \psi\n ) u}{u}{}\\
&= i\ip{\wh{T} u}{u}{} + \| \oY u\|^2 + \ip{\psi\n u}{u}{} + \ip{u}{\psi\n[0] \oY^*u}{} 
\end{align*}
So
\bgE{iT}
 i \ip{\wh{T} u}{u}{} + \| \oY^2 u\| \leqc  \|\oY^* u\|^2+ \|u\|^2 \leqc Q_{\e,\ud}^\bot(u,u) + \|u\|^2 
\enE

Now for a length one, tangential $(1,0)$ vector field $Z$,
\begin{align*}
\| Zu\|^2 &= \ip{Z^*Z u}{u}{} = \|Z^*u\|^2 + \ip{ [Z^*,Z]u}{u}{} \\
&= \|Z^*u\|^2 + \ip{i \wh{T} u}{u}{} + \ip{\psi\n u}{u}{} \\
& \leqc Q_{\e,\ud}^\bot(u,u) + \|u\|^2
\end{align*}
Thus we have establish
\bgE{10}
\| \psi\n[1,0] u \|^2 \leqc Q_{\e,\ud}^\bot(u,u) + \|u\|^2 
\enE

All that remains to be shown is that
\[ \|\oD u\|^2 \leqc Q_{\e,\ud}^\bot(u,u) +\|u\|^2.\]
But for $u\in \Cvt{\Omega}$
\begin{align*}
\| \oD u\|^2 & = \ip{\oD u}{\oD u}{} = \ip{D \oD u}{u}{} \\
&= \ip{ \oD D u}{u}{} + \ip{ [D,\oD]u}{u}{}\\
&= \ip{D u}{Du}{} +  \ip{ [D,\oD]u}{u}{}\\
& \leqc \|Du\|^2 + \ip{ \Lambda^2 u}{u}{} + \|u\|_1 \|u\|\\
& \leqc \|Du\|^2 + \|\Lambda u\|^2 + \frac{1}{\zeta} \|u\|^2 + \zeta \|u\|^2_1
\end{align*}
where $\zeta$ is chosen $0 < \zeta <<1$. Absorbing this term, yields the estimate for $\Cvt{\Omega}$. Density then implies the result.

\epf

\bgD{boxD}
We define the unbounded operators  $\Pb[\triangle,\e]$ on $L^{2,\top}(\Omega)$ by setting
 \[ \dom{\Pb[\triangle,\e]} = \left\{ u \in \Bt{\mathcal{D}} \colon \text{ there exists } f \in L^{2,\top}(\Omega) \text{ such that } Q^\bot_{\e,0} (u,u) = \ip{f}{u}{} \right\}\]
 with $\Pb[\triangle,\e] u =f$ for $u \in \dom{\Pb[\triangle,\e]}$.  We set $\Pb = \Pb[\triangle,\e]$.
\enD

From \rfL[ER]{SelfAdjoint} we can immediately deduce 
\bgL{boxD}
For $\ud>0$, $\e \geq 0$, the unbounded operators $\ud + \Pb[\triangle,\e]$ are closed, densely-defined, bijective and selfadjoint. Thus there are bounded solution operators $N^{\e,\ud}$.
\enL






\bgL{Improved}
For $0<\ud <1$ and $\e \geq 0$, if $u \in \Bt{\mathcal{D}}^\e$  with $\Pb[\triangle,\e] u =f -\ud u $ then the following estimates hold independent of $\ud$:
\bgE{I1}
 \| \Lambda^2 u\|_1 + \e\| \Lambda^2 u\|_{H^1} \leqc \|f\|_1 + \|u\|  \enE

 \bgE{I2}
  \| \Lambda^2 u\| \leqc \|f\| + \|u\|
 \enE
 \enL

\pf
First suppose that $\ud=1$, $f\in \Cct{\Omega}$ and $u \in \Cvt{\Omega}$. We note that
\[ \| \Lambda u\|^2 \leqc Q^\bot_{\e,1}(u,u)  = \ip{ f}{u}{} = \ip{ \Lambda^{-1}  f}{\Lambda u}{} \]
so
\bgE{l1} \| \Lambda u\| \leqc \| \Lambda^{-1} f \|. \enE
Now by the useful commutation property that $[\hboxb ,\Lambda^2]=0$, $[\wh{T},\Lambda^2]=\psi\n[2]$ and $[D,\Lambda^2]= \psi\n[2]$ we have
\begin{align*}
Q^\bot_{\e,1}(\Lambda^2 u, \Lambda^2 u)&= Q^\bot_{\e,1}(u, \Lambda^4 u) + \ip{ \psi\n[2] u}{ D\Lambda^2 u}{} + \ip{ Du}{\psi\n[2] \Lambda^2 u}{}\\
& \qquad + \e^2 \ip{\psi\n[2] u}{\wh{T} \Lambda^2 u}{} + \e^2 \ip{ \wh{T} u}{ \psi\n[2] \Lambda^2 u}{}\\
&=  \ip{f}{ \Lambda^4 u}{} + \ip{ \psi\n[2] u}{ D\Lambda^2 u}{} + \ip{ \psi\n[2] u}{\oD \Lambda^2 u}{}\\
& \qquad + \e^2 \ip{\psi\n[2] u}{\wh{T} \Lambda^2 u}{} 
\end{align*}
where the integration-by-parts in the last line works because if $u \in \Cvt{\Omega}$ then $\Lambda^2u \in \Cvt{\Omega}$.

From the basic estimate of \rfL{BasicEstimate} we therefore see that
\[ \| \Lambda^2 u\|_1^2 +\e^2 \|\Lambda^2 u\|_{H^1} \leqc  \|\Lambda f\| \, \|\Lambda^3 u\| + \|\Lambda^2 u\| \, \|\Lambda^2 u\|_1  + \e^2  \| \Lambda^2 u\| \, \|\Lambda^2 u\|_{H^1} \] 
and so
\[ \| \Lambda^2 u\|_1 + \e\| \Lambda u\|_{H^1}  \leqc \|\Lambda f\| + \|\Lambda^2 u\|.  \]
Now we use the observation that for $0 < \zeta <<1 $ we have
\[ \|\Lambda^2 u\| = \ip{ \Lambda^3 u}{\Lambda u}{}^{1/2}  \leqc \zeta \|\Lambda^3 u\| + \frac{1}{\zeta}  \|\Lambda u\|\]
to see that 
\[ \|\Lambda^2 u\|_1 \leqc \|\Lambda f\| + \| \Lambda u\| \leqc \| \Lambda f\|. \]
In particular this implies that
\bgE{l3}
\begin{split}
\| \Lambda^3 u \| &\leqc \|\Lambda f\|
\end{split}
\enE

Thus we can apply the interpolation result of \rfL[FN]{Interpolation} to the solution operator $N^{\e,1}$ to see that 
\[  \|\Lambda^2 u\| \leqc \|f\|.\]

Now for $0<\ud<1$ there must be  some $f\in\Cct{\Omega}$ for which any $u \in \Cvt{\Omega}$ is the $\e,\ud$-weak solution. But then $u$ is the $\e,1$-weak solution for $f+(1-\ud)u$ and so
\[ \| \Lambda^2 u \| \leqc \|f\| + |1-\ud| \, \|u\| \leqc \|f\| + \|u\|.\]
To move from \textit{a priori} estimates to genuine estimates we invoke elliptic regularity if $\e>0$ and \rfP[ER]{dReg} otherwise.

 \epf

The remaining estimates depend heavily on the following integration-by-parts computation. Suppose $Z_j$, $j=1..k+1$ are smooth vector fields that are tangent to $\partial \Omega$ at the boundary and $Z=Z_1 \dots Z_{k+1}$. We'll use $Z\n[i]$ to denote  the generic $Z$ with $k < i$. Then if $u \in \Cv{\Omega}$ is a weak solution for $f \in \Cc{\Omega}$
\bgE{Comm}
\begin{split}
 Q^\bot_{\e,\ud} (Zu,Zu) &= Q^\bot_{\e,\ud} (u,Z^*Zu) + \ip{[S\n ,Z]u}{ S\n  Zu}{} + \ip{S\n  u}{[Z^*,S\n ]Zu}{}\\
 & \qquad +\e^2 \ip{[\wh{T},Z]u}{ \wh{T} Zu}{} + \e^2 \ip{\wh{T} u}{[Z^*,\wh{T} ]Zu}{}\\
 &= \ip{f}{Z^*Zu}{} + \ip{ Z\n[k]   (S\n  +\wh{T}) u}{S\n  Z u}{} + \ip{Z\n[k]   S\n  u}{  (S\n  +\wh{T})  Z u}{} \\
 & \qquad + \e^2 \ip{ Z\n[k]   (S\n  +\wh{T} ) u}{\wh{T} Z u}{} + \e^2 \ip{ Z\n[k]   \wh{T} u}{( S\n +\wh{T}) Z u}{}  \\
 \end{split}
  \enE

The last term without the epsilon factors must be dealt with differently depending on the value of $k$. If $k=0$ then
\bgE{k0}
 \ip{ S\n u}{ (S\n  +\wh{T} ) Zu}{} = \ip{S\n  u}{S\n  Z u}{} +  \ip{ \nabla u}{S\n  Zu}{} \enE
If $k=1$ then
\bgE{k1}
 \ip{ Z\n  S\n  u}{ (S\n  +\wh{T} ) Zu}{} = \ip{ Z\n  S\n  u}{S\n  Zu}{} + \ip{ (S\n \wh{T}  + S\n[2]  ) u }{ Z\n  Zu}{}.
 \enE
If $k>1$ then
 \bgE{k2}
  \ip{ Z\n[k]   S\n  u}{ (S\n  +\wh{T} ) Zu}{} =\ip{ Z\n[k-1]   S\n  \wh{T}  u + Z\n[k-1]  S\n  u+ Z\n[k-2]   S\n[2]   u}{ Z\n  Zu}{} \enE

\[ \begin{split} &= \ip{Z\n[k]   f}{Z\n  Z u}{} + \ip{ Z\n[k]   (S\n  +\wh{T}) u}{ S\n  Z u}{} + \ip{Z\n[k]   S\n  u}{  S\n    Z u}{}\\ 
  & \qquad  +\ip{ Z\n[k-1]   S\n  \wh{T}  u + Z\n[k-1]  S\n  u+ Z\n[k-2]   S\n[2]   u}{ Z\n  Zu}{}\\
  & \qquad + \e^2 \ip{ Z\n[k]   (S\n  +\wh{T} ) u}{\wh{T} Z u}{} + \e^2 \ip{ Z\n[k]   \wh{T} u}{( S\n +\wh{T}) Z u}{}   
\end{split}\]

\bgL{Apriori}
If $0 \leq \e <<1$ and $0<\ud<1$  and $u = N^{\e,\ud} f\in \Cvt{\Omega}$ for some $f\in \Cct{\Omega}$ then we get the following a priori estimates independent of $\ud$.
\[ \|u\|_{k+2} \leqc \| f \|_k + \| u\| + \e\| u\|_{H^{k+2}}  \] 
\enL

\pf  The proof is by induction. To prove the case $k=0$, we apply \rfE{Comm} together with the relevent \rfE{k0} with each $Z=Z_1$ being a Folland-Stein vector field. Since $[ \wh{T}, Z_1] = \nabla$, we see from the basic estimate that
\[
\begin{split}
 \|Zu\|_1^2  & \leqc \|Zu\|_1 ( \|f\|+ \|u\|_{1} + \|\Lambda^2 u\| )  +\e^2 \| u\|^2_{H^2}
\end{split}\]
and so by \rfL{Improved}
\[ \|Zu\|_1 \leqc \|f\|+ \|u\| +\e \|u\|_{H^2}. \]
We now note that $\oD \rho D + D\rho \oD$ is tangent to the boundary so we can express
\bgE{Dd} D^2 =  Z\n D + \psi\n[0] \boxD + Z\n[2] + S\n +\e^2 \psi\n[0] \wh{T}^* \wh{T} \enE 
and
\bgE{oDd} \oD^2 = \oD Z\n   + \psi\n[0] \boxD + Z\n[2] + S\n + \e^2 \psi\n[0] \wh{T}^* \wh{T}.\enE
In both cases all the terms on the right hand side are already controlled. This complete the case $k=0$.

Now set $k>0$ and suppose the result is true for $j<k$. First we note for $2 \leq m < k+2$ and $u \in \Cvt{\Omega}$
 \begin{align*}
 \| \Lambda^2 u \|_{m} &\leqc \|(\ud + \Pb +\e^2 \wh{T}^* \wh{T} )  \Lambda^2 u\|_{m-2} +  \| \Lambda^2 u\| + \e \| \Lambda^2 u\|_{H^m}  \\
 &\leqc \| \Lambda^2 f\|_{m-2} + \| \psi\n[2] S\n u\|_{m-2} + \| \Pb u +\ud u\|_1 + \|u\| +\e \| u\|_{H^{m+2}}\\
 &\leqc \| \Pb u+\ud u\|_m + \|u\|_{m+1} +\e\| u\|_{H^{m+2}}
 \end{align*}

Now  we apply \rfE{Comm} together with the relevent  of \rfE{k1} or \rfE{k2} with each $Z_j$ being a Folland-Stein vector field. Since $[ \wh{T}, Z_j] = \nabla$, we see from the basic estimate that
\[
 \|Zu\|_1^2   \leqc \|Zu\|_1 ( \|f\|_{k}+ \|u\|_{k+1} + \|\Lambda^2 u\|_{k} )  +\e^2 \| u \|^2_{H^{k+2}}  
\]
and so by induction
\[\begin{split}
 \|Zu\|_1 &\leqc \|f\|_k + \|u\|_{k+1} + \|\Lambda^2u\|_k +\e \| u\|_{H^{k+2}} \\
  &\leqc   \|f\|_k + \|u\| +\e \| u\|_{H^{k+2}} . 
 \end{split}. \]
 For general derivatives we repeatedly invoke \rfE{Dd} and \rfE{oDd}.
 
 \epf

\bgL{Sobolev}

For $0 < \ud <1$,  $0 \leq \e <<1$ if $u \in \Cvt{\Omega}$ is the $\e,\ud$-weak solution for $f \in \Cct{\Omega}$, then there is the following uniform a priori estimate for all $k \geq 0$
\[ \|u\|_{H^{k+1}} \leqc \| f\|_{H^k} + \|u\|. \]
The allowable upper bound on $\e$ may depend on $k$.
\enL

\pf Again the proof is by induction. The case $k=0$ follows immediately from the basic estimate and the second part of \rfL{Improved}.

We note that from \rfL{Apriori} we get \bgE{e} \|u\|_{2} \leqc \| f\| + \|u\| + \e \|u\|_{H^{2}} \enE
with constants independent of $\e$ and $\ud$. 
Additionally from \rfL{Improved}, we see that
\bgE{LT}
\begin{split}
\| \Lambda^2 \wh{T}^m u\| &\leqc \| P \wh{T}^m u \| + \|\wh{T}^m u\| \\
& \leqc \| \wh{T}^m f\| + \| S\n T\n[m] u \| + \e^2 \| u\|_{H^{m+1}}  
\end{split}
\enE

Now since $\wh{T}^* = \wh{T} + \psi\n[0]$ we have
\begin{align*}
Q^\bot_{\e,\ud}( \wh{T}^m u,\wh{T}^m u) &= Q^\bot_{\e,\ud}( u,( \wh{T}^*)^m \wh{T}^m u) + \ip{ [\wh{T}^m ,S\n] u}{S\n \wh{T}^mu}{} \\
& \qquad  + \ip{ S\n u}{[ (\wh{T}^*)^m, S\n] \wh{T}^m u}{} \\
&\qquad + \e^2 \ip{ [\wh{T}^m ,\nabla ] u}{ \nabla \wh{T}^m u}{} + \e^2 \ip{ \nabla u}{[ (\wh{T}^*)^m, \nabla] \wh{T}^m u}{}\\
&= \ip{  f}{ (\wh{T}^*)^m \wh{T} u}{} + \ip{ (S\n T\n[m-1] +T\n[m]  )u }{S\n \wh{T}^m u}{}\\
& \qquad + \ip{S\n u}{ (T\n[m] + T\n[m-1] S\n) \wh{T}^m u}{}\\
& \qquad + \e^2  \ip{ \nabla^m u }{\nabla \wh{T}^m u}{} + \e^2 \ip{ \nabla u}{ T\n[m-1]  \nabla \wh{T}^m u}{}  \\
&= \ip{ \wh{T}^m f}{\wh{T}^m u}{} + \ip{ \nabla^m u}{S\n \wh{T}^m u}{} + \ip{ S\n u}{ T\n[m]  \wh{T}^m u}{} \\
& \qquad + \e^2 \ip{ \nabla^m u}{\nabla \wh{T}^m u}{}
\end{align*}
Now the problem term is dealt with as follows
\begin{align*}
\ip{ S\n u}{T\n[m] \wh{T}^m u}{} &= \ip{ S\n u}{ \wh{T}^m T\n[m] u + T\n[2m-1]  u }{}= \ip{S\n u}{ (\wh{T}^*)^m T\n[m] u + T\n[2m-1] u}{}\\
&= \ip{ S\n \wh{T}^m + \nabla^m u}{ T\n[m] u}{} + \ip{\nabla^m u}{\nabla^m u}{}
\end{align*}
Therefore from the basic estimate and induction we get
\[\begin{split}
 \|\wh{T}^m u \|^2_1  + \e^2 \|\wh{T}^m u\|^2_{H^1} &\leqc \big( \| f\|_{H^m} +\|u\|_{H^m} \big)\big( \|u\|_{H^m}  +\| \wh{T}^m u \|_1\big)\\ &\qquad  + \e^2  \| u\|_{H^m} \| \wh{T}^m u\|_{H^1} \end{split}\]
and so
\bgE{m}
 \| \wh{T}^m u \|_1 + \e \|\wh{T}^m u\|_{H^1} \leqc \|f\|_{H^m} + \|u\|_{H^m} 
 \enE
The case $k=1$ now follows from \rfE{m} with $m=1$, \rfE{e}, the first part of \rfL{Improved}, \rfE{LT} and the fact that
\[ \|u\|_{H^2} \leqc \| \wh{T} u \|_1 + \|u\|_2 +\| \Lambda^2 \wh{T} u\| \]
combined with the observation that for $0 \leq \e <<1$ we can absorb the last term in \rfE{e}.

We now proceed by induction and suppose $k \geq 1$ with the result being true for all $j\leq k$.  Set $P= \boxD + \ud + \e^2 \wh{T}^* \wh{T}$. If $Z$ denotes any $k$th order operator that maps $\Cvt{\Omega}$ to $\Cvt{\Omega}$, then by \rfL{Apriori}
\[ \begin{split}
\| Zu\|_2 &\leqc \| P Zu \| + \|Zu\| +\e\| u\|_{H^{k+2}} \leqc \| Zf\| + \| u\|_{H^{k+1}} \\&\leqc \|f\|_k + \|u\|+\e \|u\|_{H^{k+2}}.
\end{split} \]
We can improve this to any $k$th order differential operator by the now standard decomposition argument.

The result is then proved by noting that by  \rfE{LT} and \rfE{m},
\[ \|u\|_{H^{k+2}} \leqc \|\nabla^k u\|_2 + \| \wh{T}^{k+1} u\|_1 + \| \Lambda^2 \wh{T}^{k+1} u\|\]
using induction and absorbing the $\e \|u\|_{H^{k+2}}$ term that occurs when these terms are bounded.

\epf

\bgT{Main}
For $0 <\ud<1$, for $f \in L^{2,\top}(\Omega)$ there is a unique $u \in \dom{\Pb}$ such that $\Pb u = f - \ud u$. Furthermore if $f \in S^k$ (respectively $H^k$) then $u \in S^{k+2}$ (respectively $H^{k+1}$) and there are the following estimates independent of $\ud$
\begin{align*}
\| u \|_{k+2}  &\leqc \|f \|_k + \|u\| \\
\| u\|_{H^{k+1}} &\leqc \|f\|_{H^k} + \|u\| 
\end{align*}
 \enT
 
 \pf
 Hypoellipticity of the solution operators follows from  \rfP[ER]{dReg} and \rfL{Sobolev}. This also implies the Sobolev estimates. Once we have hypoellipticity, the Folland-Stein estimates follow immediately from the \textit{a priori} estimates of \rfL{Apriori}.
 
 \epf 
 
 \bgC{Main}
 For $\ud =0$, if $\Pb u \in S^k$ then $u \in S^{k+2}$  and there is a uniform estimate
 \[ \|u\|_{k+2} \leqc \| \Pb u\| + \|u\|. \]
  \enC

\bgR{1}
The upper bound $\ud<1$ was chosen fairly arbitrarily. The reason for chosing a bound is to ensure uniformity across choice of $\ud$. However, for any $\ud>0$ the above arguments can be used to show that $\ud+\Pb$ is hypoelliptic and satisfies the estimates of \rfT{Main}.
\enR

\section{Tangential regularity and estimates} \setS{APt}

The goal of this section is to prove \text{a priori} estimates for the operators $\ud+\Pt$ with $\ud \geq 0$. Unfortunately, these operators are not globally subelliptic and the the possible presence of cohomology on the leaves plays havoc on the delicate estimates that were obtained in earlier work. In fact, to prove the estimates for $\Pt$ we shall first have to establish hypoellipticity for the operators $ \ud +\Pt $ with $0<\ud<1$. 

Once more, throughout this section we shall suppose that $(\Omega,w,\eta)$ satisfies (A1), (A3) and (B).

\bgD{DQ} We define the sesquilinear forms
\begin{align*}
 Q^{\top}_{\ud,\e}(u,v) &= \ip{\hdb u}{\hdb v}{} + \ip{\hdbs u}{\hdbs v}{} + \ip{\oD u}{\oD v}{} \\
 & \qquad + \ud \ip{u}{v}{} + \e^2 \ip{\nabla u}{\nabla u}{} 
 \end{align*}
 with $Q^\top = Q^\top_{0,0}$. Then set $\Bt{\mathcal{E}}^\e$ to be the closure of $\Cct{\Omega}$ in the $Q^\top_{1,\e}$-norm.
 \enD
 
 
 
  
\bgD{Pt}
We define
\[ \begin{split} \dom{\Pt[\triangle,\e]} = \Big\{ u\in \Bt{\mathcal{E}}^\e \colon &\text{ there exists }f\in L^{2,\top}(\Omega) \text{ such that } \\ &Q^\top_{\e,0}(u,v) = \ip{f}{v}{} \text{ for all }v\in \mathcal{E} \Big\} \end{split}\]
and set $\Pt[\triangle,\e] u = f$. We set $\Pt = \Pt[\triangle,0]$.
\enD

From \rfL[ER]{Properties} we can immediately deduce that
\bgL{Properties}
For $\ud>0$ the operator $\ud+\Pt[\triangle,\e]$ is densely-defined, closed, bijective and selfadjoint as unbounded operator on $L^{2,\top}(\Omega)$. Thus there exists a bounded solution operator $N^{\e,\ud}$.
\enL 





 \bgL{apriori1}
Suppose $\ud>0$, $\e \geq 0$ and  $f$ and $u= N^{\e,\ud} f$ are both in $\Cct{\Omega}$. Then
\[ \| \Lambda^2 u\| + \|\Lambda \oD u\| + \e \|\Lambda u\|_{H^1} \leqc \|f\| \]
with the constant independent of $\e$. 
\enL

\pf
We clearly have the basic estimate that 
\[ \| \Lambda u \|^2 + \|\oD u\|^2 + \e^2 \|u\|^2_{H^1} \leqc \frac{1}{\ud} Q^\top_{\ud,\e}(u,u) \leqc \frac{1}{\ud} \|f\| \|u\| \leq \frac{1}{\ud^2} \|f\|^2.\]

Repeating the arguments of \rfL[TV]{Improved} establishes that
\bgE{L2}
 \|\Lambda^2 u \| \leqc \|f\| + \|u\|. 
 \enE
 with constant independent of $\ud$.  Now we  integrate by parts argument with any vector field $X \in \Psi_1$.
 \begin{align*}
 Q^\top_{\ud,\e}(Xu,Xu) &= \ip{f}{X^*Xu}{} + \ip{ [\db, X]u}{\db Xu}{} + \ip{u}{ [X,\db] X u}{} \\
 & \qquad  + \ip{ [\dbs, X]u}{\dbs Xu}{} + \ip{u}{ [X,\dbs] X u}{} \\ &\qquad + \ip{ [\oD,X]u}{\oD Xu}{} +\ip{ \oD u}{[X, \oD] X u}{} \\
& \qquad + \e^2 \ip{ [\nabla, X]u}{ \nabla X u}{} + \e^2 \ip{\nabla u}{[\nabla,X] X u}{}\\
& \leqc[\ud] \|f\| \|\Lambda^2 u\| +  \|\Lambda^2 u\|^2 +\|\Lambda u\| \, \|\Lambda \oD u\| + \|\Lambda u\|^2 + \e^2 \|\Lambda u\|\, \|Xu\|_{H^1} \\
& \leqc[\ud]   \left(  \|\Lambda^2 u \| + \|\Lambda \oD u\| + \e \|X  u\|_{H^1} \right) \|f\|.
 \end{align*}
Applying the basic estimate and summing over a spanning set of such $X$, yields the result.

\epf

\bgC{SimpleReg}
$\dom{\Pt} \subset \wh{S}^2$.
\enC

\pf
Apply the elliptic regularization argument of \rfP[ER]{dReg}.

\epf

Unfortunately, \rfL{apriori1} does not give a coercive basic estimate in this case. Counter-examples to the existence of such an estimate were constructed in \cite{hladky1}. However, we make the key observation that formally $\db \boxb  = \boxb \db $. The implication is that after using the the decomposition of \rfC[CP2]{decomposition}, we can essentially then use the regularity and estimates of \rfS{TV} to control certain derivatives of the solutions. 

\bgL{Ddom}
Suppose that $u \in \dom{\Pt}$ and $\Pt u \in \dom{\oD}$, then $\oD u \in \dom{\Pb}$.
\enL

\pf
From the definition of $\Pt$, we see that there is some $f \in L^{2,\top}(\Omega)$ such that
\[Q^\top_{0,0}(u,v)= \ip{f}{v}{}\]
for all $v \in \Cct{\Omega}$. Now suppose instead that $v \in \Cvt{\Omega}$. Then by repeated integration by parts
\[ Q^\bot_{0,0}(\oD u,v) = Q^\top_{0,0}(u, D v) + \ip{\psi\n u}{\psi\n v}{}.\]
Now by \rfC{SimpleReg}, we see that $u \in \wh{S}^2$ for all $X \in \Psi_2$. So we can integrate by parts one more time to see
\[ Q^\bot_{0,0}(\oD u, v) = \ip{f}{Dv}{} + \ip{\psi\n[2] u}{v}{} = \ip{\oD f + \psi\n[2] u}{v}{}.\]
But this implies that $\oD u \in \dom{\Pb}$.
 
\epf

 \bgL{Dbar}
Suppose $u\in \Cct{\Omega} \cap \dom{\Pt}$ and $\ud \geq 0$,  then 
\[ \| \oD u\|_1 \leqc \|\Pt[\ud] u\| + \|\Lambda u\| \] 
and for $k>0$
\[ \| \oD u \|_{k+1} \leqc \|\Pt[\ud] u \|_k + \| \Lambda^2 u\|_{k-1}. \]
\enL

\pf We'll prove it for $\ud=0$, the case $\ud>0$ is almost identical.

From \rfL[CP2]{Properties}  we see that, applied to smooth forms, $\oD \Pt  - \Pb \oD = \psi\n[2]$. Furthermore if $u \in \dom{\Pt}$ then $\oD u \in \dom{\Pb}$, but this implies that by  \rfT[TV]{Main}
\[ \|\oD u\|_{j+2} \leq \|\Pb \oD u\|_j + \|\oD u\|.\]
Therefore for $k>0$
\[ \|\oD u\|_{k+1} \leq \| \oD \Pt u \|_{k-1} + \|\Lambda^2 u\|_{k-1}. \]
Now for $k=0$, we note that
\begin{align*}
 \|\oD u\|_1^2 &\leq Q^\bot(\oD u , \oD u) + \|\oD u\|^2\\
 &= \ip{ \hboxb \oD u}{\oD u}{} + \ip{\oD^* \oD u}{\oD^* \oD u}{} + \|\oD u\|^2\\
 &= \ip{\hboxb u}{\oD^* \oD u}{} + \ip{[\hboxb,\oD]u}{\oD u}{}  + \ip{\oD^* \oD u}{\oD^* \oD u}{} + \|\oD u\|^2\\
 & = \ip{\Pt u}{\oD^* \oD u}{} + \ip{\psi\n[2] u}{\oD u}{} +\|\oD u\|^2\\
 & \leqc \| \Pt u\| \|\oD u\|_1 + \|\Lambda u\|\, \|\oD u\|_1 + \|\oD u\|^2
 \end{align*}
 so
 \[ \|\oD u\|_1 \leqc \|\Pt u\| + \|\Lambda u\|. \]
\epf

\bgL{1hypoelliptic}
For $\ud>0$, the operator $\ud+\Pt$ is hypoelliptic and satisfies the estimate
\[ \|\Lambda^2 u\|_k + \|\oD u\|_{k+1}  \leqc[\ud] \|\ud u+ \Pt u\|_k \]
for all $u\in \Cct{\Omega}$.
\enL

\pf
Suppose that $f \in \Cct{\Omega}$ and $\Pt u =f -\ud u$. First we note that $\Pb \oD = \oD \Pt + \psi\n[2]$. Now by the \text{a priori} estimates of \rfL{apriori1} and the elliptic regularization argument of \rfP[ER]{dReg}, we see that $\Lambda^3 u, \Lambda^2 \oD u \in L^{2,\top}$ and hence $ \oD u \in \dom{\Pb}$ and $(\ud+\Pb) \oD u = \oD f - \psi\n[2] u\in L^{2,\top}$ But this implies that $\oD u \in S^2$, which in turn implies that $\Lambda_\ud^2 u \in S^1$. This then implies $(\ud+\Pb) \oD u \in S^1$ and so $\oD u \in S^3$. Thus $\Lambda_\ud^2 u \in S^2$ and we can continue to boot strap our way up to see that $u \in S^\infty \subset \Cct{\Omega}$.

The estimates then follow from \rfL{Dbar} and the observation that $\Lambda^2 \leqc[\ud] \Lambda_\ud^2 \leqc[\ud] \Lambda^2$.

\epf

Having established the hypoellipticity of $\ud+\Pt$ for $\ud>0$, we now attend to the task of finding \text{a priori} estimates for $\Pt$ itself. The possibility of cohomology on the leaves, makes this substantially more difficult as an estimate for $\| \hboxb u\|_k$ does not immediately imply an estimate for $\|u\|_k$. To overcome this issue, we shall again use the idea of interpolation on the foliating leaves. Our key observation is the following. For for $u \in \Cct{\Omega}$
\bgE{goal} \ip{ D^k u}{D^k u}{}  = \ip{\Lambda D^k u}{\Lambda^{-1} D^k u}{} \leq \e \| \Lambda D^k u\| + \frac{1}{\e} \|\Lambda^{-1} D^k u\| \enE
The goal is then to establish estimates for the terms on the right hand side of \rfE{goal}.

\bgL{Minus1} For all $u\in \Cct{\Omega}$,
\[ \| \Lambda^{-1} Du \| \leqc \| \Pt[F] u \| + \| \Lambda^2 u \| + \|\oD u\|. \]
where $\Pt[F]$ is the formal operator $\Pt[F] = \hboxb + D\oD$.
\enL

\pf  First note that 
\bgE{D}
\| \Lambda^{-1} D u\|  \leq \sup\limits_{ \varphi \in \Coi[\Omega]} \dfrac{ \left| \ip{ Du}{\varphi}{} \right|}{\| \Lambda \varphi \|} =  \sup\limits_{ \varphi \in \Coi[\Omega]} \dfrac{ \left| \ip{u}{\oD \varphi}{} \right|} {\| \Lambda \varphi\|} 
\enE
Now since $1+\Pt$ is hypoelliptic, if we set $N=N^{0,1}$, then $ N \varphi \in \Cc{\Omega} \cap \dom{\Pt}$. But then $\oD N \varphi \in \dom{\oD^*}$ and so $\oD N \varphi \in \Cvt{\Omega}$. Thus
\begin{align*}
 \ip{u}{\oD \varphi}{} &= \ip{u}{ (1+ \Pt ) N \oD \varphi}{} \\
 &= \ip{ \hboxb u}{N \oD \varphi}{} + \ip{\oD u}{\oD N \oD \varphi}{} + \ip{u}{N \oD \varphi}{} \\
 &= \ip{\hboxb u}{\oD N \varphi}{} + \ip{\hboxb u}{[N ,\oD] \varphi}{} \\ & \qquad + \ip{\oD u}{\oD^2 N \varphi}{} + \ip{\oD u}{\oD [N,\oD] \varphi}{}\\ & \qquad + \ip{u}{\oD N \varphi}{} + \ip{u}{[N,\oD]\varphi}{}\\
  &= \ip{\hboxb u}{\oD N \varphi}{} + \ip{\hboxb u}{N [\oD,1+\Pt]N \varphi}{} \\ & \qquad + \ip{D \oD u}{\oD N \varphi}{} + \ip{\oD u}{ \oD N [\oD,1+\Pt] N \varphi}{} \\ &\qquad + \ip{u}{\oD N \varphi}{} + \ip{u}{N [\oD,1+\Pt]N \varphi}{}
\end{align*}
Now 
\[ [\oD,\Pt] = \psi\n[2] + [\oD, D] \oD = \psi\n[2] + (S\n + \psi\n[2]) \oD.\]
Note that $\Lambda^2$ is a well defined, differential operator which acts in directs tangent to the foliation, $[ \oD,\Lambda^2]=\psi\n[2]$ and  $[\hdb,\Lambda^2]=0=[\hdbs,\Lambda^2]$. Thus we can once again apply a commutation argument to the basic estimate to see that for $\ua \in \Cc{\Omega} \cap \dom{\Pt}$
\begin{align*}
\| \Lambda^3 \ua \|^2 +\| \oD \Lambda^2  \ua \|^2 & \leqc Q^\top_{0,1}(\Lambda^2 \ua, \Lambda^2 \ua)  \\
& = Q^\top_{0,1} (\ua, \Lambda^4 \ua ) + \ip{\psi\n[2] \ua}{\oD \Lambda^2 \ua}{} 
\end{align*}
and so  there is an estimate
\[ \| \Lambda^2 \oD \ua \| \leqc \| \Lambda (1+\Pt) \ua \|. \]
Applying this to $\ua = N\varphi$ we have the estimate 
\[ \| [\oD,\Pt] N \varphi  \|  \leqc \| \Lambda \varphi \|. \]
However we also have 
\[ \| \oD N \ua \| \leq \| \ua \|, \qquad \|N \ua\| \leq \|\ua\| \]
for all smooth $\ua$,  and so
\[ \left| \ip{Du}{\varphi}{}  \right| \leqc \left( \| \Pt[F] u\| + \|\Lambda^2 u\| + \|\oD u\| \right) \| \Lambda \varphi \| \]
and the result is proved.

\epf

\bgL{Useful} 
The following commutation and integration properties hold:
\begin{itemize}
\item $[ \oE^k, \psi\n] = \sum\limits_{j<k} \psi\n \oE^{j} = \sum\limits_{j<k} \oE^j \psi\n$.
\item $[ \oE^k, \nabla] = \sum\limits_{j<k} \nabla \oE^{j} = \sum\limits_{j<k} \oE^j \nabla$.
\item $\oE^* = \oE + \psi\n[0]$
\end{itemize}
\enL

\bgT{MainAP}
Suppose $0 < \ud <<1$ and $u \in \Cct{\Omega} \cap \dom{\Pt}$ such that $\Pt u =f -\ud u $. Then we have the a priori estimate uniform over $\ud$
\[ \|\Lambda^2 u\|_k + \|\oD u\|_{k+1}  \leqc \|f\|_k + \|u\| \]
for all $k \geq 0$.
\enT

\pf The proof will follow a contorted induction argument. The case $k=0$ follows from \rfL{Dbar} and  \rfL{apriori1} combined with the observation that $\| \Pt u\| \leq \|\Pt[1]u\| + \|u\|$. Now suppose the result is true for all $0\leq j<k$. 

First we note that by \rfL{Dbar} and the inductive hypothesis, we have
\[ \| \oD u \|_{k+1} \leqc \|  f \|_{k} + \|\Lambda^2 u\|_{k-1} \leqc \|f\|_k + \|u\|. \]
Additionally, we note that
\[ \| \hboxb u\|_k = \| f - D\oD u -\ud u \|_k \leqc   \|f\|_k +\ud \|u\|_k+ \|u\|\]
so
\[ \|\Lambda^2 u\|_k \leqc \|f\|_k + \|u\|_k. \]
By the inductive hypothesis, the only derivative we now need to control to establish the estimate is $\|D^k u\| \leqc \|f\|_k + \|u\|$. This is the meat of the argument. Now
\begin{align*}
Q^\top(\oE^k u,\oE^k u) &= \ip{ \hdb \oE^k u}{ \hdb \oE^k u}{} + \ip{\hdbs \oE^k u}{\hdbs \oE^k u}{} + \ip{ \oD \oE^k u}{\oD \oE^k u}{} \\
& \qquad + \ud \ip{ \oE^k u }{\oE^k u}{} \\
&= Q^\top( u,(\oE^k)^* \oE^k u) + \ip{ \sum\limits_{j <k} \psi\n \oE^j  u}{ \psi\n \oE^k u}{}  \\ &\qquad + \ip{ \nabla \sum\limits_{j<k} \oE^ju }{\oD \oE^k u}{}  + \ip{\oD u}{ \sum_{j<k} \oE^j \nabla \oE u}{}\\
&= \ip{ \oE^k f}{\oE^k u}{} +\ip{\psi\n \sum\limits_{j<k} \oE^j u}{\psi\n \oE^k u}{} + \ip{\nabla \sum\limits_{j<k} \oE^j u }{\oD \oE^k u}{} \\ &\qquad +\ip{ \sum\limits_{j \leq k} (\oE^*)^j \oD u}{\nabla \oE^{k-1}  u}{} + \ip{\sum\limits_{j<k} (\oE^*)^j \oD u}{ \nabla \sum\limits_{j<k} \oE^j u}{}\\
&= \ip{ \oE^k f}{\oE^k u}{} +\ip{\psi\n S\n[k-1] u}{\psi\n \oE^ku}{} + \ip{\nabla S\n[k-1] u }{\oD \oE^k u}{} \\ &\qquad +\ip{\oD  \sum\limits_{j \leq k} \oE^j u + \nabla S\n[k-1] u }{\nabla \oE^{k-1}  u}{} + \ip{\nabla S\n[k-1] u}{ \nabla S\n[k-1] u}{}
\end{align*}
So by the inductive hypothesis
\begin{align*}
 \left| Q^\top(\oE^k u,\oE^k u) \right| & \leqc \left( \|\oE^k f\| + \| S\n[k] u\| +\|\Lambda^2 S\n[k-1] u\|\right)\left( \|\Lambda \oE^k u\| + \|\oD \oE^k u\| \right)\\ &\qquad  + \| u\|_k^2 + \|\Lambda^2 S\n[k-1] u\|^2\\
 & \leqc   \left( \| f\|_k + \| u\|_k \right)\left( \|\Lambda \oE^k u\| + \|\oD \oE^k u\| \right) + \| u\|_k^2 + \|f\|^2_{k-1}
\end{align*}
Hence by the basic estimate
\[ \|\Lambda \oE^k u\| + \|\oD \oE^k u\| \leqc \|f\|_k + \|u\|_k. \]
Now we are almost done for
\[ \Lambda D^k \leqc \Lambda \oE^k + \text{controlled terms} \]
so
\[ \|\Lambda D^k u\| \leqc \|f\|_k + \|u\|_k. \]
From \rfL{Minus1} and the inductive hypothesis we get
\begin{align*}
 \|\Lambda^{-1} D^k u\| &\leqc \|\Pt[F] D^{k-1} u \| + \|\Lambda^2 D^{k-1} u \| + \|\oD D^{k-1} u\| \\ 
 &\leqc \|f\|_{k-1} + \ud \| u \|_{k-1} + \| [\Pt[F],D^{k-1} ]u\| + \|u\| + \| [\oD,D^{k-1}] u\|\\
 &\leqc \|f\|_{k-1} + \|\Lambda^2 S\n[k-1]u\| +\|u\|  +\| \Lambda^2 S\n[k-2] u\| \\
 &\leqc \|f\|_{k-1} + \|u\|
 \end{align*}
Using our leaf interpolation method, we see
\begin{align*}
\|D^k u\|^2 &\leq \e \|\Lambda D^k u\|^2 + \frac{1}{\e} \| \Lambda^{-1} D^k u\|^2 \\
& \leqc \e \left( \|f\|_k^2 + \|u\|_k^2 \right) + \frac{1}{\e} \left( \|f \|_{k-1}^2 +\|u\|^2 \right)
\end{align*}
Thus 
\[ \|u\|_k^2 \leq \e \|u\|^2_k + \|f\|^2_k + \|u\|^2. \]
Choosing $\e$ sufficiently small and absorbing the $\|u\|_k$ term on the right, then yields the result.

\epf

\bgT{Main}
Suppose that $0 < \ud <<1$, $u \in \dom{\Pt}$ and $\ud u + \Pt u \in S^k$ for $k\geq 0$. Then we have the following regularity results: $u \in S^k$, $\Lambda^2 u \in S^k$, $\oD u \in S^{k+1}$, $\varrho u \in S^{k+1}$ and $\varrho^2 u \in S^{k+2}$. Furthermore there is a uniform estimate independent of $\ud$ of the form
\[ \| \Lambda^2 u\|_k + \|\varrho u \|_{k+1} + \|\varrho^2 u\|_{k+2} + \|\oD u \|_{k+1} \leqc \| \Pt u \|_k + \|u\|. \]
Furthermore if we additionally insist that $u \bot \Ker{\Pt}$ then these results extend to the case $\ud=0$.
\enT

\pf
Most of this theorem has already been proved. Since we know the operators $\ud + \Pt$ are surjective and hypoelliptic, it easily follows that the \textit{a priori} estimates of \rfT{MainAP} are genuine estimates. It just remains to show the additional weighted regularity. Again from hypoellipticity and surjectivity it suffices to show the estimate holds for smooth $u \in \dom{\Pt}$.

Now if $u \in \Cct{\Omega}$, then $\varrho u$ vanishes on the boundary. Alternatively  phrased, $\varrho u$ satisfies Dirichlet boundary conditions. We can then repeat the arguments from \rfS{TV}, for the operator $\ud+\Pt$ to see that 
\[ \| \varrho u \|_{k+1} \leqc \| \Pt (\varrho u) \|_{k-1} + \| \varrho u\| \]
where we compute $\Pt (\varrho u)$ formally and if $k=0$ replace $k-1$ by $0$. But
\[ \| \Pt \varrho u \|_{k-1} \leqc \| \varrho \Pt u \|_{k-1} + \| u\|_k \leqc \|\Pt u\| + \|u\|.\]
The same argument can be applied to $\varrho^2 u$ to obtain
\[ \|\varrho^2 u\|_{k+2} \leqc \| \varrho^2 \Pt u \|_k + \| \rho u \|_{k+1} + \| u \|_k .\]
Thus the result holds for $0 < \ud <<1$. The case for $\ud=0$ follows from \rfL[ER]{0Reg}.

\epf

\section{Regularity and Estimates for \boxb} \setS{MN}

Throughout this section we shall suppose that $u$ and $f$ are $(0,q)$-forms with $1 \leq q \leq n-2$  and $\ud >0$. Now we set $\mathcal{D} = \Cc{\Omega} \cap \dom{\db} \cap \dom{\dbs}$

 \bgL{density} If $\Omega$ has noncharacteristic boundary and satisfies (A1) and (A2) then the set $\mathcal{D}$ is dense in $\dom{\db} \cap \dom{\dbs}$ in the graph norm, \[u \mapsto \|u\| + \| \db u\| + \| \dbs u\|.\]
\enL

Since $\Omega$ has noncharacteristic boundary and satisifies  (A1) and (A2), the density result can be proved using a method almost identical to Lemma 4.3.2 \cite{Shaw:B}. The equivalent notion to splitting a form into complex tangent and normal pieces is just the decomposition $u = u^\top + \theta^\bt{0} \wedge u^\bot$. The details are lengthy but standard.

\bgD{Q}
The sesquilinear forms $Q^{\ud,\e}$ are defined by
\[ Q^{\ud,\e}(u,v) = \ip{\db u}{\db v}{} + \ip{\dbs u}{\dbs v}{} + \ud\ip{u}{v}{} + \e \ip{\nabla u }{\nabla v}{}.\]
The spaces $\Bt{\mathcal{D}}^\e$ are now the closures of $\mathcal{D}$ under that $Q^{1,\e}$-norm.
\enD

\bgL{AP} Under conditions (A1),  (A3) and (B),
for $u \in \mathcal{D}$ and $\ud,\e \geq 0$ such that $Q^{\ud,\e}(u,v) = \ip{f}{v}{}$ for all $v \in \Bt{\mathcal{D}}^\e$ 
\[ \|\Lambda u\| \leqc \|f\| +\|u\| \]
\enL

\pf
Now for $u \in \mathcal{D} $
\begin{align*}
\|u\|^2+ \| \Lambda u\|^2 & \leqc \|\Lambda u^\top\|^2 +\| \Lambda u^\bot\|^2 \leqc Q_{0,0}^\top(u^\top,u^\top) + Q_{0,0}^\bot(u^\bot,u^\bot) + \|u\|^2 \\
&=Q^{0,0}(u,u) + \ip{\psi\n u^\top}{u^\bot}{} + \ip{\psi\n u^\bot}{u^\top}{} + \|u\|^2\\
& \leqc Q^{\ud,\e}(u,u) + \|\Lambda u\| \, \|u\| + \|u\|^2\\
&= \| f\| \, \|u\| + \|\Lambda u\| \, \|u\| + \|u\|^2\\
& \leqc \left( \|\Lambda u\| + \|u\| \right)  \left( \|f\| + \|u\| \right)
\end{align*}

\epf

\bgC{1reg}
Under conditions (A1), (A3) and (B), $\dom{\boxb} \subset \wh{S}^1$
\enC

\pf
First note that if $u \in \dom{\boxb}$ then $u \in \dom{ \ud + \boxb}$ for any $\ud \geq 0$. Now apply \rfP[ER]{dReg} with $\mathcal{X} = \wh{S}^1$, $\mathcal{Y}= L^2(\Omega)$ and $P=\boxb$  to see that \[u = N^\ud (\ud  +\boxb )u  \in \wh{S}^1.\]

\epf

\bgL{DiagSA} Under conditions (A1), (A3) and (B), 
the operator 
\[ \boxD u = \Pt u^\top + \theta^\bt{0} \wedge \Pb u^\bot\]
with $\dom{\boxD} = \left\{ u \colon u^\top \in \dom{\Pt}, u^\bot \in \dom{\Pb} \right\}$ is self-adjoint.
\enL

\pf
This follows immediately from \rfL[TV]{boxD} and \rfL[APt]{Properties}.

\epf

\bgL{ED} Under conditions (A1), (A3) and (B) $\dom{\boxb} = \dom{\boxD}$.
\enL

\pf
\rfT[TV]{Main} and \rfT[APt]{Main} imply that $\dom{\boxD} \subset \dom{\boxb}$. To prove the reverse inclusion we see that for $u \in \wh{S}^1 \cap \dom{\boxb}$ and $v \in \dom{\boxD}$
\begin{align*}
\ip{u}{\boxD v}{} &= \ip{u}{\boxb v +\psi\n v}{} = \ip{\boxb u + \psi\n u}{v}{}.
\end{align*}
Since $\boxD$ is self-adjoint, this implies that $u \in \dom{\boxD^*} = \dom{\boxD}$ and so we have
\[ \dom{\boxb} \cap \wh{S}^1 \subset \dom{\boxD} \subset \dom{\boxb}.\]
But by \rfC{1reg}, $\dom{\boxb} \subset \wh{S}^1$.

\epf

For a pair $(\Omega,w)$ such that (A1) holds we denote by $K$ the points  $p \in \Omega$ such that $H_w \bt{w}$ is not constant on the leaf of $w$ through $p$.  If (A2) holds then $K$ is contained in the interior of $\Omega$ and we can define
\[  \mathscr{B}^\infty = \left \{ \xi  =\xi(w,\bt{w}) \in \Cc{\Omega} \colon 1-\xi \in \Coi[\Omega], \overline{K} \subset \{ \xi =0 \} \right\} .\]
Thus $\mathscr{B}^\infty$ is the collection of smooth functions $\xi$ depending only on $w,\bt{w}$ such that $\xi =1$ on a neighborhood of $\partial \Omega$ and (A3) holds on the support of $\xi$.
 
\bgD{k2}
The spaces $S^{k;j}$ are defined as follows: for $\xi \in \mathscr{B}^\infty$ 
\[\begin{split}
S^{k;j} = \Big\{ u \in S^k \colon & \varrho^j u \in S^{k+j} \text{for $j=1,\dots,k$}, (\xi u)^\bot \in S^{k+j}, \\ &  \oD (\xi u^\top) \in S^{k+j-1}, \Lambda^j (\xi u) \in S^{k} \Big\}
\end{split}
 \]
 with associated norm
 \[ \| u\|_{k;j} = \sum\limits_{m=0}^j \|\varrho^m u\|_{k+m} + \|\xi u^\bot\|_{k+j} + \| \oD (\xi u^\top)\|_{k+j-1} + \|\Lambda^j (\xi u)\|_k.\]
\enD
Up to equivalence of norms, this definition is independent of the choice of $\xi$.

 \bgL{APboxb}
 Under conditions (A1), (A3) and (B), for $u \in \dom{\boxb} \cap \Cc{\Omega}$ there is a uniform $\text{a priori}$ estimate
 \[ \|u\|_{k;2} \leq \|\boxb u \|_k + \|u\|. \]
 \enL
 
 \pf From \rfT[TV]{Main} and \rfT[APt]{Main} we see
 \begin{align*}
 \|u\|_{k;2}&  \leqc \| u^\top\|_{k;2} + \|\theta^\bt{0} \wedge u^\bot\|_{k+2} \leqc \| \boxD^\top u^\top \|_k + \| \boxD^\bot u^\bot\|_k + \|u\|\\
 & \leqc \|\boxD^\top u^\top\|_k  + \|\boxb u\|_k +\|\psi\n u^\top\|_k + \|u\|\\
 & \leqc \|\boxD u^\top\|_k  + \|\boxb (\theta^\bt{0} \wedge u^\bot)\|_k + \|u\|\\
&\leqc \|\boxb u\|_k + \|\psi\n u^\bot\|_k    + \|u\|\\
& \leqc \|\boxb u\|_k  + \|\boxD^\bot u^\bot\|_{k-1} + \|u\| 
\end{align*}
By repeating this argument, we see that
\begin{align*}
\|u\|_{k;2}& \leqc \|\boxb u\|_k +\|\boxD^\bot u^\bot\| +\|u\|\\
& \leqc \|\boxb u\|_k + \|\psi\n u^\bot\| + \|u\|\\
& \leqc \|\boxb u\|_k + \|u\|
\end{align*}
where the last inequality follows from \rfL{AP}.

 \epf
 
 \bgL{Regd} Under conditions (A1), (A3) and (B),
 if $\boxb u +\ud u \in S^k$ then $u \in S^{k;2}$.
 \enL
 
 \pf
 The proof is by induction again. First we note that if $u\in \dom{\boxb}$ then $u \in \dom{\boxD}$ by \rfL{ED}. Then \rfT[TV]{Main} and \rfT[APt]{Main} imply that $u \in S^{0;2}$. Thus the result is true for $k=0$. 

Suppose the result is true for $j<k$. If $u \in \dom{\boxb}$ and $\boxb u +\ud u\in S^k$ then $u \in S^{k-1;2}$ and
\[ (\ud + \boxD)u  = \ud u+ \boxb u + \theta^\bt{0} \wedge (\psi\n  u^\top) + \psi\n u^\bot  \]
Decomposing  we get
\[  (\ud + \boxD^\top)u^\top  \in S^{k}, \qquad (\ud +\boxD^\bot)u^\bot \in S^{k-1}.\]
The second of these implies that $u^\bot \in S^{k+1}$. This implies that we can improve the first to $(\ud+\boxD^\top)u^\top \in S^k$. This implies that $u^\top \in S^{k;2}$. Thus we can improve the second to $(\ud+\boxD^\bot)u^\bot \in S^k$ which implies that $u^\bot \in S^{k+2}$. This establishes the result.
 
 \epf
 
 \bgC{Hypoellipticity}
 For $\ud>0$, the operators $\ud + \boxb$ are hypoelliptic.
 \enC

 We summarize our results so far in the following theorem:

\bgT{SubMain}
Under conditions (A1), (A3) and (B), if $1 \leq q \leq n-2$ and $0 < \ud << 1$ then for any $(0,q)$-form $f \in L^2(\Omega)$  and $\ud>0$ there is a unique $u \in \dom{\boxb}$ such that $ \boxb u=f- \ud u$. Furthermore if $f \in S^k(\Omega)$ then $u \in S^{k;2}(\Omega)$ and there is a uniform estimate independent of $\ud$
\[ \|u\|_{k;2} \leqc  \|f\|_k +\|u\| \]
\enT

 \pf
 The operators $\ud + \boxb$ are self-adjoint, injective and have closed range. Thus for $f \in L^2(\Omega)$ there is a unique $u \in \dom{\boxb}$ with $\boxb u =f -\ud u$. The regularity of $u$ follows from \rfL{Regd}. From \rfL{APboxb} and the observation that $\boxb$ is continuous from $S^{k;2}$ to $S^k$ we see that
 \[ \| u\|_{k;2} \leqc \|f\|_k +\ud\|u\|_k + \|u\|. \]
For sufficiently small $\ud$ we can therefore absorb the $\ud \|u\|_k$ term.

\epf

\bgT{Main}
Suppose $\Omega$ is a noncharacteristic smoothly bounded domain such that the triple $(\Omega,\theta,w)$ satisfies (A1) and (A2). If $1 \leq q \leq n-2$ and $\boxb$ has closed range in $L^2_{0,q}(\Omega)$ then for any $(0,q)$-form $f \in L^2(\Omega)$ such that $f \bot \Ker{\boxb}$ there is a unique $u \bot \Ker{\boxb}$ such that $\boxb u=f$. Furthermore if $f \in S^k(\Omega)$ then $u \in S^{k;2}(\Omega)$ and there is a uniform estimate
\[ \|u\|_{k;2} \leqc  \|f\|_k. \]
\enT

\pf The content here is that we must relax the assumptions in \rfT{SubMain}  from (A3) to the weaker (A1) and allow (B) to fail in the interior of $\Omega$.  However $C(w)= \{ p \colon dw \wedge d\bt{w} =0 \} $ is closed so $K=C(w) \cap \Bt{\Omega}$ is a compact set which by condition (A2) lies inside $\Omega$.  Therefore we can construct a nest of smoothly bounded open sets $U_j$ such that 
\[ K \subset U_1 \subset \subset U_2 \subset\subset  U_3 \subset \subset \Omega .\]
Let $\Omega_j = \Omega \backslash U_j$ for $j=1,2,3$. Since (A2) holds near the boundary, we can choose these sets so that (A2) holds on each $\Omega_j$. Thus the triples $(\Omega_j,\theta,w)$ all satisfy (A), (A3) and (B).

As we are assuming that $\boxb$ has closed range and \boxb is self-adjoint, there is a decompostion
\[ L^2(\Omega) = \rng{\boxb} \oplus \Ker{\boxb}.\]
Thus if $f \bot \Ker{\boxb}$ then there is a unique $u \bot \Ker{\boxb}$ such that $\boxb u = f$.

Now suppose additionally that $f \in S^k(\Omega)$ and choose  $\xi \in \mathscr{B}$  such that $\xi=1$ on $\Omega_3$ and $\xi=0$ on $U_1$. For $\ud>0$ there is a unique $u_\ud \in \dom{\boxb}$ such that
\[ \boxb u_\ud = f-\ud u_\ud.\]
By standard interior regularity results (\cite{Folland:H}, \cite{Shaw:B}), we see that $u_\ud \in S^{k+2}(U_3)$ and 
\bgE{Interior}
\|u_\ud\|_{k+2,U_3} \leqc \|f\|_k + \|u_\ud\|.
\enE
Now, it can easily  be seen that $\xi u_\ud \in \dom{\boxb[{\Omega_1}]}$  and
\[ \boxb (\xi u_\ud) = \xi  f -\ud \xi u_\ud + \zeta \ns u \]
where $\zeta$ is a smooth real-valued function that vanishes on both $\Omega_3$ and $U_1$. Therefore by \rfT{SubMain} we see
\begin{align*}
 \| \xi u_\ud\|_{k;2} &\leqc \| \xi f +\zeta \ns u\|_k  + \| u_\ud \| \\
 & \leqc \|f\|_k + \|u_\ud\|_{k+1,U_3}  \\
 & \leqc \|f\|_k + \|u_\ud\|
 \end{align*}
 Combining with \rfE{Interior}, we see that each $u_\ud \in S^{k;2}$ and there is a uniform estimate independent of $\ud$ that
 \[ \|u_\ud \|_{k;2} \leqc \|f\|_k + \|u_\ud \|.\]
 We can then apply the regularity result of \rfL[ER]{0Reg} together with the closed range assumption to see that $u \in S^{k;2}$ and that the desired estimate holds.

\epf

 

\section{Proofs of the theorems}\setS{PI}

The following result from functional analysis will prove very useful in allowing us to move results between degrees and pseudohermitian structures.

\bgL{ClosedRange}
The following are equivalent.
\renewcommand{\theenumi}{\arabic{enumi}}
\begin{enumerate}
\item $\boxb$ has closed range in $L^2_{(p,q)}$.
\item $\db$ and $\dbs$ have closed range in $L^2_{(p,q)}$. 
\item $\db$ has closed range in $L^2_{(p,q)}$ and $L^2_{(p,q+1)}$.
\item $\db$ has closed range in $L^2_{(p,q+1)}$ and $\dbs$ has closed range in $L^2_{(p,q-1)}$. 
\end{enumerate}
\enL

\pf From standard results in functional analysis on the closed range properties of adjoint operators it follows that (2), (3) and (4) are equivalent.

Now suppose (1) that $\boxb$ has closed range in $L^2_{(p,q)}$. Then
\[ L_{(p,q)}^2(\Omega) = \rng{\boxb} \oplus \Ker{\boxb}.\]
Now $\Ker{\boxb} \subset \Ker{\dbs} = \overline{\rng{\db}}^\bot$. Thus 
\[ \overline{\rng{\db}}\subset \rng{\boxb} \subset \rng{\db} \oplus \rng{\dbs.} \]
But this clearly implies that $\rng{\db}$ is closed in $L^2_{(p,q)}$. A virtually identical argument implies to $\rng{\dbs}$.

Next suppose that (2) and hence (3) and (4) hold, so that both $\db$ and $\dbs$ have closed range on the relevant bidegree forms. Then 
\[ L_{(p,q)}^2(\Omega) = \rng{\db} \oplus \Ker{\dbs} = \rng{\db} \oplus \rng{\dbs} \oplus H \]
 where $H$ is the orthogonal complement of $\rng{\dbs}$ in $\Ker{\dbs}$. Since standard results imply that $\rng{\dbs}^\bot =\Ker{\db}$ we see that $H \subset \Ker{\db} \cap \Ker{\dbs} = \Ker{\boxb}$.
 Thus we can decompose \[u= u_1 + u_2 + u_H\] with $u_1 \in \rng{\db}$, $u_2 \in \rng{\dbs}$ and $u_H \in H$.  
 Now suppose $u \in \dom{\db} \cap \dom{\dbs}$ and $\ u \bot \Ker{\boxb}$. Then $u_H=0$, $u_1 \in \dom{\dbs} \cap \Ker{\db}$ and $u_2 \in \dom{\db}\cap \Ker{\dbs}$. But since both $\db$ and $\dbs$ have closed range we see
 \[ \|u\|^2 \leq \|u_1\|^2 + \|u_2\|^2 \leqc \|\dbs u_1\|^2 + \|\db u_2\|^2. \]
 But for $u \in \dom{\boxb }$ the far right hand side is equal to $\ip{\boxb u}{u}{}$. Thus $\boxb$ has closed range in $L^2_{(p,q)}$.
 
\epf

It is frequently desirable to work with a particular fixed pseudohermitian structure rather than a rescaled form. Thus it is useful to see how things change under a rescaling. While the regularity results for the Kohn Laplacian itself do not easily hold up under this rescaling, as we shall see the closed range property and results for related inhomogenous $\db$ equation to rescale well. 

The  key observation is as follows: for $(0,q)$-forms $\varphi$, $\psi$
\begin{align*}
 \ip{\mu_x \varphi}{\mu_x \psi}{x} &= \int_\Omega \aip{ \mu_x \varphi}{ \mu_x \psi}{x} \eta^\au \wedge (d \eta^\au)^n \\
 &= \int_\Omega e^{-qx} \aip{ \varphi}{\mu}{x}  e^{(n+1)x} \eta \wedge d\eta^n = \ip{e^{(n+1-q)x} \varphi}{\psi}{}  
 \end{align*}
 where $\mu_x$ is defined back in \rfE[NL]{mu}. From this it is easy to see that
\[ \db (\mu_x \varphi)  = \mu_x (\db \varphi), \qquad \dbs (\mu_x \varphi)  = e^{-x}  \mu_x \left(  \dbs \varphi   -(n+1-q) \overline{\db x} \vee \varphi \right). \]

\bgC{ClosedRange}
If $\boxb$ has closed range on $(p,q)$-forms then for all $x \in \Cc{\Omega}$ with $x >0$, $\boxb[x]$ has closed range also.
\enC

\pf
It is sufficient to show that $\db$ has closed range for the pseudohermitian form $\eta^\au$ in $L^2_{(p,q)}$ and $L^2_{(p,q+1)}$. Suppose that $\db (\mu_x u_n) \to f$. Then $\db u_n \to \mu_{-x} f$. But since $\db$ has closed range for $\eta$ itself, there is some $u \in L^2(\Omega)$ such that $\db u = \mu_{-x} f$. But then $\db \mu_x u =f$ and so $f$ is in the range of $\db \mu_x$.

\epf

\bgC{Kernel}
If $\boxb$ has closed range, then the dimension of $\Ker{\boxb[x]}$ is independent of choice of smooth $x>0$.
\enC

\pf
Since the assumption implies that $\Ker{\db} =\rng{\db} \oplus \Ker{\boxb}$ the result follows from the observation that $\mu_x$ preserves both $\Ker{\db}$ and $\rng{\db}$.

\epf

We now have all the ingredients to prove the theorems from the introduction.

\vspace{10pt}

\noindent \textbf{Proof of Theorem A:}

Theorem A is just \rfT[MN]{Main}.

\epf



\noindent \textbf{Proof of Theorem B:}

Normalize $\theta = \eta^\au$ so that $(\Omega,\theta,w)$ additionally satisfies (D) near the boundary. From \rfC{ClosedRange}, this normalization  preserves the closed range condition. 

Let $N$ be the Neumann operator for $\boxb[0,q]$ with this rescaled pseudohermitian form. Then $v=N \mu_x f$ solves $\boxb[0,q] v = \mu_x f$ uniquely with $v \bot \Ker{\boxb[0,q]}$. Set $u = \dbs v$. Then since $v \in \dom{\boxb[0,q]}$ and $0 =\db \boxb[0,q] v= \boxb[0,q+1] \db v $ it can easily be see that \[ \db v \in \Ker{\boxb[0,q+1]} \subset \Ker{\db} \cap \Ker{\dbs} .\] Therefore 
\[ \db u = \db \dbs v = \db \dbs v +\dbs \db v = \boxb[0,q+1] v = \mu_x f.\] In particular, we now have that \[ \db \mu_{-x} u =f. \]
For the estimates, we first observe that the decomposition of $\dbs$ in \rfL[CP2]{db} implies that $\dbs$ maps $S^{k;2} \cap \dom{\dbs}$ continuously into $S^{k;1}$. The estimates then follow easily from Theorem \ref{A} and the observation that each $\mu_{x}$ is an isomorphism on all function spaces that depend only on derivatives of component functions.

\epf

\noindent \textbf{Proof of Theorem C:}




As the Kohn Laplacian is globally subelliptic on compact \spc manifolds the lack of cohomology implies the existence of a strictly positive smallest eigenvalue on each leaf. Continuity of eigenvalues and the lack of cohomology on each boundary leaf implies the existence of smoothly bounded set $\partial \Omega \subset \Omega_b \subset  \Bt{\Omega}$ foliated by level sets of $w$ such that $\Omega_b$ is open in $\overline{\Omega}$ and there is a global, strictly positive lower bounded on eigenvalues for leaves contained within $\Omega_b$. Since $dw \wedge d\bar{w} \ne 0$ on $\partial \Omega$ we can normalize the pseudohermitian form so that (A3) holds  on $\Omega_b$.

This all implies  the existence of an estimate of the form
\bgE{L2} \| \Lambda u^\top \|_{L^{2,\top}(\Omega_b)} \leqc \ip{ \hboxb u^\top}{u^\top}{L^{2,\top}(\Omega_b)}.\enE

Therefore we can apply the basic estimates of \rfS{TV} and \rfS{APt} to the domain $\Omega_b$ and combine with \rfE{L2} to see that 
\bgE{Estb}
\begin{split}
 \| \Lambda u^\top \|^2 + \| u^\bot\|_1^2 &\leqc \ip{\Pt u^\top}{ u^\top}{} + \ip{\Pb u^\bot}{ u^\bot}{} \\
 & \leqc \ip{ \boxb  u}{  u}{} + \ip{\psi\n  u^\top}{u^\bot}{} + \ip{\psi\n  u^\bot}{ u^\top}{} \\
 & \leqc \ip{\boxb  u}{ u}{} + \| u^\bot\| \, \|\Lambda  u^\top \|
 \end{split}\enE
 for $u \in \dom{\boxb}$ on $\Omega_b$. 
 
 Now let $\xi_1 \subset \xi_2$ be smooth bump functions compactly supported within $\Omega$ whose supports contain $\Bt{\Omega} - \Omega_b$ and set $\zeta = 1- \xi_1$. 
 
Applying \rfE{Estb} to $\zeta u$ then implies that there is an estimate on $\Omega$ itself of the form
 \bgE{Est2}  \|  \zeta u^\top \|+ \|\zeta u^\bot\|_1 \leqc \ip{\boxb u}{u}{} + \| \xi_2 u\|_1 + \| \zeta u^\bot \|\enE
 for all $u \in \dom{\boxb}$ on $\Omega$.

 Now suppose that $u_n$ is a sequence in $\dom{\boxb}$ such that $u_n \perp \Ker{\boxb}$, $\boxb u_n \to 0$ and $\|u_n\|=1$. From \rfE{Est2}, this implies that $\zeta u^\bot_n$ is a bounded sequence in $S^1$. Since $\Bt{\Omega}$ is compact, by applying the Rellich lemma and passing to a subsequence we can assume that $ \zeta u_n^\top$ converges in $L^2$.  A similar argument using the interior estimate
 \[  \| \xi_2 u\|_1 \leqc \| \boxb u\| + \|u\| \]
implies that a subsequence of $\xi_2 u_n$ converges. The estimate \rfE{Est2} then implies that an appropriate subsequence of $\zeta u_n^\top$ is $L^2$-Cauchy and so converges.  Now since $\{ \zeta =1\} \cup \{ \xi_2 =1 \} \supset \Bt{\Omega}$ we can deduce that after passing to a subsequence $u_n$ itself is Cauchy and so converges in $L^2$.  Thus $u_n \to u$ for some $u \in L^2$ and $\boxb u_n \to 0$. Since $\boxb$ is a closed operator, this implies that $u \in \Ker{\boxb}$. But this is a contradiction, and hence we must have an estimate
\[ \| u\| \leqc \|\boxb u\|, \qquad \text{for all }u \perp \Ker{\boxb}.\]
This estimate is equivalent to $\boxb$ having closed range. 
 A similar argument this time assuming that $u_n \in \Ker{\boxb}$ shows that every bounded sequence in $\Ker{\boxb}$ has a convergent subsequence. Thus the $L^2$ unit ball in $\Ker{\boxb}$ is compact as so $\Ker{\boxb}$ is finite dimensional.
 
 Thus for this normalized pseudohermitian form $\boxb$ is Fredholm on $L_{0,q}^2(\Omega)$. From \rfC{ClosedRange} and \rfC{Kernel}, we see that $\boxb$ for the original pseudohermitian form is also Fredholm.




Next we show that $\boxb$ is hypoelliptic up to the boundary when (A3) holds. From \rfT[MN]{Main} it suffices to show that if $\boxb u =0$ then $u \in \Cc{\Omega}$. Standard interior estimates imply that $u \in C^\infty(\Omega)$. It remains to show that $u$ is smooth up to the boundary. As the boundary leaves have no cohomology we can find a smooth neighborhood $U$ of the boundary such that $\Omega_b = U \cap \Omega$ satisfies (A1), (A3), (B) and none of the foliating leaves in $\Omega_b$ have cohomology. The closed range result implies that $\Omega_b$ also satisfies (E).  Let $\xi$ be a smooth positive function supported in $U$ that  is identically equal to $1$ on $\partial \Omega$. Since $u \in \Ker{\boxb}$, it follows that $f = \boxb (\xi u)  \in \Cc{\Omega_b}$ and $\xi u \in \dom{\boxb}$ on $\Omega_b$.

Now from \rfC[CP2]{decomposition} we see that there is a first order operator $L \in \Psi_1$ such that
\[ \boxb v = \Pt v^\top + L v^\bot + \theta^\bt{0} \wedge ( \Pb v^\bot + L^* v^\top ) \]
for all $v \in \dom{\boxb}$ on $\Omega_b$.  Thus $\xi u$ solves the following system
\begin{align*}
 \Pt v^\top  &= f^\top  - L v^\bot  \\
 \Pb v^\bot &= f^\bot - L^* v^\top.
\end{align*}
As $\dom{\boxb} \subset \wh{S}^1$ and vitally that $\Pt$ is injective on $\Omega_b$, we can now run an induction argument using \rfC[TV]{Main} and \rfT[APt]{Main} to see that $\xi u \in \Cc{\Omega}$. For if $L^*(\xi u) \in S^{k-1}$ then \rfC[TV]{Main} implies that $(\xi u)^\bot \in S^{k+1}$. Then \rfT[APt]{Main} implies that both $(\xi u) \in S^k$ and $L^*(\xi u)^\top \in S^k$. 

\epf







 \noindent \textbf{Proof of Theorem D:}
 
 For $p=0$, this follows immediately from \rfT[ID]{noco} and Theorem C. To move to the case $p>0$,  we recall that $dz^1$, \dots, $dz^{n+1}$ yields a global holomorphic trivialization of $\Lambda^{1,0}(M)$. Since $\boxb$ has closed range on $(0,q)$-forms implies by \rfL{ClosedRange} that $\db$ has closed range in $(0,q)$ and $(0,q+1)$-forms. The presence of global holomorphic trivialization immediately implies that $\db$ has closed range in $(p,1)$ and $(p,2)$. For if $ \db (dz^I \wedge u_n ) \to \phi$ then $\phi$ can be written as $dz^I \wedge u$ and $u_n \to u$. But this implies $u= \db v$ for some $v$ and hence $\phi= \db (dz^I \wedge v)$. But then \rfL{ClosedRange} can be applied again to see that $\boxb$ has closed range on $(p,q)$-forms. 
 
 That $\boxb$ has finite kernel follows a very similar argument. Finite dimensional kernel on $(0,q)$-forms is equivalent to  $\rng{\db}$ having a finite dimensional complement in $\Ker{\db}$. The global holomorphic trivialization immediately extends this later condition to $p>0$.

\epf

\noindent \textbf{Proof of Theorem E:}

From Theorem D we see that $\boxb$ is Fredholm and from Theorem C that it is hypoelliptic on $\overline{\Omega}$. A theorem due to Shaw \cite{Shaw:P} states that the extra assumptions that $w=z_1$ and $\varrho$ is strictly convex imply that the Kohn Laplacian is actually injective.  Theorem B then implies smooth solvability of $\db$ on $(0,q)$-forms. For the case $p>0$ we again use the global holomorphic trivialization.

\epf

\appendix
\section{Functional Analysis and Elliptic Regularization} \setS{ER}

In this section, we prove the elliptic regularization results that we shall need in various places throughout the paper.

Here are the standing assumptions we shall make throughout the section

\begin{itemize}
\item $\mathcal{D}$ is a linear subspace of $\Cc{\Omega}$ such that
\[ \Coi \subset \mathcal{D} \subset \Cc{\Omega}\]
\item $Q(u,v)$ is a first order, symmetric, positive sesquilinear form on $\mathcal{D}$
\item $\mathcal{D}$ is the closure of $\mathcal{D}$ in the $1+Q$-norm.
\item $\mathcal{X}$ and $\mathcal{Y}$ are Hilbert spaces such that $\Cc{\Omega} \subset \mathcal{X} \subset \mathcal{Y}  \subset L^2(\Omega)$ and \[\|u\| \leq \|u\|_\mathcal{Y} \leq \|u\|_\mathcal{X}\] 
\item $\Cc{\Omega}$ is dense in $\mathcal{Y}$.
\end{itemize}

We define \[ Q^{\ud,\e}(u,v) =  Q(u,v) + \ud\ip{u}{v}{} + \e \ip{\nabla u}{\nabla v}{} \]
on $\Bt{\mathcal{D}}^\e$, the closure of $\mathcal{D}$ in the $Q^{1,\e}$-norm. From this we define second operator operators $P_{\e,\ud}$ by 
\[\begin{split} \dom{P_{\ud,\e} } = \Big\{ u \in \Bt{\mathcal{D}}^\e \colon &\text{ there exists $f\in L^2(\Omega)$ such that } \\ & Q^{\e,\ud}(u,v) = \ip{f}{v}{} \text{ for all $v \in \Bt{\mathcal{D}}^\e$} \Big\} \end{split}\]
with $P_{\ud,\e} u =f$. We set $P_\ud = P_{0,\ud}$, $P = P_{0,0}$ and $\Bt{\mathcal{D}} = \Bt{\mathcal{D}}^0$.



\bgL{SelfAdjoint}
An unbounded, closed, densely defined, symmetric operator, $F$, on a Hilbert space $\mathcal{X}$ that is bijective from $\dom{F} \to \mathcal{X}$ is selfadjoint.
\enL

\pf Since $F$ is surjective, it has closed range. Therefore since $F$ is injective, closed and densely defined, we have the bound
\[ \|u \| \leqc \|F u\| \]
for all $u \in \dom{F}$. This implies that $F^{-1}$ is a bounded operator $\mathcal{X} \to \mathcal{X}$. But $F^{-1}$ is also symmetric and any bounded, symmetric operator is selfadjoint. But this implies that $F$ itself is selfadjoint  (\cite{RieszNagy}, p.312).
 
\epf

\bgL{Properties}
For $\ud>0$, $P_{\e,\ud}$ is a closed, densely defined, bijective, selfadjoint unbounded operator on  $L^2(\Omega)$.
\enL

\pf
The operator $P_{\e,\ud}$ is densely-defined as $\Coi[\Omega] \subset \dom{P_{\e,\ud}}$. From the positivity of $Q$ we immediately see
\[\| u\|^2 \leq \frac{1}{\ud} Q_{\ud,\e}^\top(u,u).\] To show that it's closed, suppose that $u_n \in \dom{P_{\e,\ud}}$ is a sequence such that $u_n \to u$ and $P_{\e,\ud} u_n \to f$ in $L^2(\Omega)$. But  the estimate then implies $u_n$ converges in $\Bt{\mathcal{D}}^\e$ in the $Q^{\e,\ud}$-norm. Thus $u \in \Bt{\mathcal{D}}$ and
\[ Q^{\e,\ud}(u,v) = \lim Q^{\e,\ud}(u_n,v) = \lim \ip{P_{\e,\ud} u_n}{v}{} = \ip{f}{v}.\]

The estimate also implies both that $P_{\e,\ud}$ is injective and using the Riesz representation theorem applied to $v \mapsto \ip{f}{v}{}$ on $\Bt{\mathcal{D}}^\e$, that there exists a weak solution $u \in  \Bt{\mathcal{D}}^0$ to $Q^{\e,\ud}(u,v) = \ip{f}{v}{}$ for all $f \in L^{2}(\Omega)$. Therefore we clearly have that the operators are bijective. As $Q$ is symmetric so is $P_{\e,\ud}$. Self-adjointness then follows from \rfL{SelfAdjoint}.

\epf

Since $P_{\e,\ud}$  is bijective for $\ud>0$, we can always construct a bounded inverse operator $N^{\e,\ud}$. Again we use $N^\ud =N^{0,\ud}$. 

We'll base our main regularity result on the following key lemma:

\bgL{WC}
Suppose $u_n$ is a bounded sequence in $\mathcal{X}$ that converges weakly in $L^2$ to some $u$, then $u \in \mathcal{X}$ and 
\[ \| u\|_{\mathcal{X}} \leq \lim \inf \|u_n\|_{\mathcal{X}}  \] 
\enL

\pf
Since $\mathcal{X}$ is a Hilbert space, it is reflexive and so $\|u_n\|$ bounded  in $\mathcal{X}$ implies that $u_n$ has a subsequence $u_{n_k}$ that converges weakly to some $\ub \in \mathcal{X}$ and
\[ \| \ub\|_{\mathcal{X}} \leq   \lim \inf \|u_n\|_{\mathcal{X}}. \]

Now for $v \in L^2(\Omega)$,
\[ \ip{u}{v}{} \leq \| u \| \,\|v\| \leq \|u\|_{\mathcal{X}} \|v\| ,\qquad \text{for all $u \in \mathcal{X}$.}\]
Therefore the map $u \mapsto \ip{u}{v}{}$ is a bounded linear functional on $\mathcal{X}$. But this implies that $u_{n_k} \to \ub$ weakly in $L^2(\Omega)$ and so $\ub =u$.

\epf

Next we prove our main elliptic regularization theorem. 

\bgP{dReg}
Suppose $\ud>0$ and the following \text{a priori} estimate holds uniformly in $\e$ : if $u \in \mathcal{D}$ and $Q^{\ud,\e}(u,v)=\ip{f}{v}{}$ for all $v \in \Bt{\mathcal{D}}^\e$ then
\[  \| u\|_{\mathcal{X}}  \leqc \|f\|_{\mathcal{Y}} + \|u\|. \]
Then if $f \in \mathcal{Y}$ then $N^\ud f \in \mathcal{X}$ and
\[  \|N^\ud f\|_\mathcal{X} \leqc \| f\|_\mathcal{Y} + \frac{1}{\ud} \|Nf\| \]
\enP

\pf 
For $f\in L^2(\Omega)$
\[  \ud \| N^{\ud,\e} f \|^2 \leq Q^{\ud,\e}(N^{\ud,\e}f , N^{\ud,\e}f )  \leq  \ip{ f}{N^{\ud,\e}f }{} \leq  \|f\| \|N^{\ud,\e}f \|\]
so
\[ \| N^{\ud,\e}f \| \leq  \frac{1}{\ud} \|f\|, \qquad  \text{for all $f \in L^2$, $\e \geq 0$.} \]
Now for $\e >0$ there is an elliptic estimate
\[ \| u\|^2_{H^1} \leq \frac{1}{ \min (\e,\ud) } Q^{\ud,\e}(u,u), \qquad \text{for all $u\in \Bt{\mathcal{D}}^\e$.} \]
Thus we can apply elliptic regularity to see that if $f \in \Cc{\Omega}$ we can deduce that $N^{\ud,\e}f \in \Cc{\Omega}$. 
Applying the \text{a priori} estimate we see that for $\e>0$ and $f\in \Cc{\Omega}$
\bgE{Key}
 \| N^{\ud,\e} f \|_{\mathcal{X}} \leq C_k  \left( \|f\|_\mathcal{Y} + \|N^{\ud,\e} f\| \right) \leq  C\|f\|_\mathcal{Y} + \frac{C}{\ud} \|f\| 
 \enE
with $C$ independent of $\e,\ud,f$.

Once again suppose that $f \in \Cc{\Omega}$. As $\e \to 0$, the sequence $N^{\ud,\e} f$ is bounded in $\mathcal{X}$. But for $g \in \mathcal{D}$
\[ \left| Q^\ud (N^{\ud,\e} f - N^\ud f, g) \right| = \e \left| \ip{\nabla N^{\ud,\e}f }{\nabla g}{} \right| \leq \e \| N^{\ud,\e} f\|_{H^1} \|g\|_{H^1} \]
Now 
\begin{align*}
\|N^{\ud,\e} f\|_{H^1} &\leq \sqrt{  \frac{1}{ \min (\e,\ud) } Q^{\ud,\e}(N^{\ud,\e} f, N^{\ud,\e} f) }  \leq \sqrt{  \frac{1}{ \min (\e,\ud) } }  \sqrt{ \| f\| \, \|N^{\ud,e} f\| }\\
& \leq \frac{1}{ \sqrt{\ud} \sqrt{ \min (\ud,\e) } } \|f\|
\end{align*}

Thus
\[ \left| Q^\ud (N^{\ud,\e} f - N^\ud f, g) \right| \leq  \frac{\e}{\sqrt{\ud} \sqrt{ \min (\e,\ud) }} \|f\| \|g\|_{H^1} \longrightarrow 0  \qquad \text{as $\e \to 0$.}\]
This is sufficient to show  that \bgE{Almost} Q^\ud(N^{\ud,\e}f -N^\ud f, g) \to 0\enE for all $ g \in \mathcal{D}$.

However
\begin{align*}
Q^\ud(N^{\ud,\e} f-N^\ud f ,&N^{\ud,\e} f - N^\ud f) = \e  \ip{\nabla N^{\ud,\e}f }{\nabla N^{\ud,\e} f}{} - Q^\ud(N^{\ud,\e}f,N^\ud f)\\
& \qquad  \qquad \qquad \qquad  + Q^\ud(N^\ud f, N^\ud f)\\
& \leq  \frac{\e}{\ud \min (\e,\ud) } \|f\|^2  + \|f\| \, \|N^{\ud,\e}f \| + Q^\ud(N^\ud f,N^\ud f)\\
& \leq \frac{1}{\ud} \|f\|^2 + Q^\ud(N^\ud f,N^\ud f)
\end{align*}
so $N^{\ud,\e} f -N^\ud f$ is bounded in $\Bt{\mathcal{D}}^0$ as $\e \to 0$. This combined with \rfE{Almost} implies that $N^{\ud,\e}f -N^\ud f$ converges weakly to 0 in $\Bt{\mathcal{D}}^0$.

Thus $N^{\ud,\e} f$ is a bounded sequence in $\mathcal{X}$ that converges weakly in $L^2(\Omega)$ to $N^\ud f$. Hence by \rfL{WC}, we see that $N^\ud f \in \mathcal{X}$ and
\bgE{Key2} \| N^\ud f\|_\mathcal{X}   \leqc \|f\|_\mathcal{Y} +\frac{1}{\ud} \|N^\ud f\|. 
\enE

Now suppose $f\in \mathcal{Y}$, then by the density assumption there is a sequence of smooth forms $f_n \to f$ in $\mathcal{Y}$. Now as $n \to \infty$, $N^{\ud} f_n$ is therefore a bounded sequence in $\mathcal{X}$ that converges to $N^{\ud}f$ in $L^2(\Omega)$. By \rfL{WC} this implies that $N^{\ud}f \in \mathcal{X}$  and \rfE{Key2} holds for all $f \in \mathcal{Y}$.

\epf

For $\ud=0$, we no longer have surjectivity for $P$. If we assume that $P$ has closed range however, we can often extend the preceding result partially to $\ud=0$. It only proves regularity for the solutions to $Pu=f$ that are orthogonal to $\Ker{P}$. In the case of the Kohn Laplacian, this all we expect. Examples of  forms in the kernel of $\boxb$ that are not even in $S^1(\Omega)$ were constructed in \cite{hladky1}.

 \bgL{0Reg}
Suppose we assume in addition to the standing assumptions that
\begin{itemize}
\item $P$ has closed range
\item The Neumann operators $N^\ud$ are hypoelliptic
\item  There is a an \text{a priori} estimate for $u \in \mathcal{D} \cap \dom{P}$
\[ \|u\|_{\mathcal{X}} \leqc \|Pu\|_{\mathcal{Y}} + \|u\|. \]
\item $P$ is continuous as a map from $\dom{P} \cap \mathcal{X}$ to $\mathcal{Y}$
\end{itemize}
Then if $Pu=f \in \mathcal{Y}$ and $u \bot \Ker{P}$ then $u \in \mathcal{X}$ and
\[ \|u\|_{\mathcal{X}} \leqc \|f\|_\mathcal{Y} . \]
 \enL
 
 \pf Since we are assuming $P$ has closed range there is a constant $c_0$ such that
\[ \|u\| \leq c_0 \|P u\|, \qquad \text{for all $u \in \dom{P}$ with $u \bot \Ker{P}$}\] and therefore there is a Neumann operator $N \colon \rng{P} \to (\Ker{P})^\bot \cap \dom{P}$ with 
\begin{align*}
 N P &= \one, \qquad \text{ on $(\Ker{P})^\bot \cap \dom{P}$}\\
 P N &= \one, \qquad \text{ on $\rng{P}$.}
 \end{align*}
 Furthermore we have the estimate
 \[ \|Nf \| \leq c_0 \|f\|, \qquad \text{for all $f\in \rng{P}$.} \]
Now note that for $f \in L^2(\Omega)$, $N^\ud f \in \dom{P_\ud} = \dom{P}$. Furthermore
\[ P N^\ud f = (P_\ud -\ud)N^\ud f = f - \ud N^\ud f. \]
Therefore  $f \in \rng{P}$ implies that $N^\ud f \in \rng{P} = (\Ker{P})^\bot.$
In particular, we now have
\[ \| N^\ud f \| \leq c_0 \left( \|f\| + \ud \|N^\ud f\| \right)\]
and so for small $0<\ud<<1$
\[ \| N^\ud f\| \leqc  \|f\|. \]

Suppose $f \in \Cc{\Omega}$ then by assumption $N^\ud f \in \Cc{\Omega}$. Now by the \text{ a priori} estimate we see
\bgE{key3}
\begin{split}
 \| N^\ud f\|_\mathcal{X} &\leqc  \left(  \| P N^\ud f\|_\mathcal{Y} + \|N^\ud f\| \right) \\
 & \leqc  \| f\|_\mathcal{Y} + \ud\|N^\ud f\|_\mathcal{Y} +  \|N^\ud f\| .
\end{split} \enE
so by absorbing the middle term for $0 < \ud <<1$ we get
\bgE{key4}
\|N^\ud f\|_\mathcal{X} \leqc \|f\|_\mathcal{Y} +\|N^\ud f\|
\enE
 Now for $f \in \mathcal{Y}$, there is a sequence $f_n \in \Cc{\Omega}$ such that $f_n \to f$ in $\mathcal{Y}$. As $n\to \infty$, the sequence $N^\ud f_n$ is bounded in $\mathcal{X}$. But since $\ud \|N^\ud g\| \leq \|g\|$ for all $g \in L^2(\Omega)$, the sequence  converges strongly in $L^2(\Omega)$ to $N^\ud f$. By \rfL{WC}, this shows that $N^\ud f \in \mathcal{X}$ and \rfE{key4} holds as a genuine estimate.
  
Now  for $f \in \rng{P} \cap \mathcal{Y}$
 \[ \left| Q (N^\ud f -Nf,g)\right| = \ud \left| \ip{N^\ud f}{g}{} \right| \leq \ud \| N^\ud f\| \, \|g\| \leqc  \|f\| \, \|g\| \to 0 \]
 for all $g \in \Bt{\mathcal{D}}^0$. Since both $N^\ud f$ and $N f$ are in $(\Ker{P})^\bot$ this implies that $N^\ud f \to Nf$ strongly in $L^2(\Omega)$. This implies that $N^\ud f$ is a bounded sequence in $\mathcal{X}$ that converges weakly to $Nf$ in $L^2(\Omega)$. By \rfL{WC} this implies that $Nf \in \mathcal{X}$ and
 \[ \|Nf\|_\mathcal{X} \leq \lim \inf \|N^\ud f\|_\mathcal{X} \leqc \|f\|_\mathcal{Y}. \]

 \epf


\bibliographystyle{plain}
\bibliography{References}

\end{document}